\documentclass{amsart}

\usepackage{amssymb}

\numberwithin{equation}{section}

\newtheorem{theo}{Theorem}[section]
\newtheorem{defi}[theo]{Definition}
\newtheorem{lemma}[theo]{Lemma}

\newtheorem{corol}[theo]{Corollary}
\theoremstyle{remark}
\newtheorem{remark}[theo]{Remark}

\newcommand{\mc}[1]{\mathcal{#1}}
\newcommand{\mbb}[1]{\mathbb{#1}}
\newcommand{\uq}{\mc{U}_{q}}
\newcommand{\uqt}{\mc{U}_{\tilde q}}
\newcommand{\qt}{\tilde q}

\newcommand{\phl}{\pi_\lambda}
\newcommand{\phlb}{\pi_{\bar \lambda}}
\newcommand{\obx}{\hat{x}}

\newcommand{\htr}{H^{\lambda}_{\tau,\rho}}
\newcommand{\fus}{F^{\lambda}_{\upsilon,\sigma}}

\newcommand{\tc}[1]{\widetilde{\mathcal{#1}}}
\newcommand{\hg}{\hat \gamma}

\title[Ruijsenaars' hypergeometric function and the modular double]
{Ruijsenaars'
hypergeometric function and
the modular double of $\uq(\mathfrak{sl}_2(\mbb{C}))$}
\author{Fokko J. van de Bult}
\address{Korteweg-de Vries Institute for Mathematics, University of Amsterdam, Plantage Muidergracht 24, 1018 TV Amsterdam, The Netherlands}
\email{fjvdbult@science.uva.nl}
\date{\today}

\begin{document}

\begin{abstract}
Simultaneous eigenfunctions of two Askey-Wilson second order difference 
operators
are constructed as formal matrix coefficients of the principal
series representation of the modular double of the quantized
universal enveloping algebra $\uq(\mathfrak{sl}_2(\mbb{C}))$. 
These eigenfunctions are
shown to be equal to Ruij\-senaars' hypergeometric 
function under a proper parameter correspondence.
\end{abstract}

\maketitle

\section{Introduction}
The main goal of this article is to construct a solution of two commuting 
 Askey-Wilson second order difference equations 
 using representation theory of the modular double of the
quantized universal enveloping algebra
 $\uq$ of $\mathfrak{sl}_2(\mathbb{C})$. 
Furthermore we relate this solution
to Ruijsenaars' hypergeometric function from \cite{Ruijs2}.

By \cite{Masuda} there exist three inequivalent $*$-structures on
$\uq$, one associated to the
real form $\mathfrak{su}(2)$ of $\mathfrak{sl}_2(\mathbb{C})$, 
one associated to $\mathfrak{su}(1,1)$, and
one to $\mathfrak{sl}_2(\mbb{R})$.
Koornwinder \cite{Koorn}, Noumi \& Mimachi \cite{Noumi}, and Koelink 
\cite{someone} have shown that the
Askey-Wilson polynomials arise as matrix coefficients of $*$-unitary irreducible representations of
 $\uq(\mathfrak{su}(2))$. To prove these results they used the 
fact that the
 Askey-Wilson second order difference operator arises as the radial part
of the quantum Casimir in $\uq$
 when calculated with respect to Koornwinder's \cite{Koorn}
twisted primitive elements. In \cite{SenK} and \cite{Stokman} 
Koelink and Stokman constructed  
the trigonometric Askey-Wilson functions as matrix coefficients of $*$-unitary irreducible representations of 
$\uq(\mathfrak{su}(1,1))$. In this paper we consider matrix coefficients
of $\uq(\mathfrak{sl}_2(\mathbb{R}))$-representations.

An essential tool is the embedding of 
$\uq(\mathfrak{sl}_2(\mbb{R}))$ in Faddeev's \cite{Faddeev}
 modular double of $\uq$.
The modular double consists of two commuting copies of the quantized
universal enveloping algebra of $\mathfrak{sl}_2(\mbb{C})$ with 
deformation parameters $q=e^{\pi i w_1/w_2}$ and $\qt=e^{\pi i w_2/w_1}$ 
($w_1,w_2 \in \mbb{R}_{>0}$), respectively.
Kharchev et al.\ \cite{Kharchev} made the crucial observation that
the algebraic version $\phl$ of the principal series representation of 
$\uq(\mathfrak{sl}_2(\mbb{R}))$
on the space $\mc{M}$ of 
meromorphic functions on $\mbb{C}$ can be extended to a representation
of the modular double on the same space. In the same article they construct 
generalized Whittaker functions as matrix coefficients of $\phl$.

We construct joint eigenvectors to the action under $\phl$ of two
 commuting twisted primitive elements 
(one for each copy of the quantized universal enveloping algebra of 
$\mathfrak{sl}_2(\mbb{R})$ inside the modular double) in terms of 
Ruijsenaars' \cite{Ruijs1} hyperbolic gamma function. 
The action of the two commuting quantum Casimir elements in the modular double
shows that the corresponding matrix coefficients, for which we have an 
explicit integral representation, satisfy Askey-Wilson second order 
difference equations in step directions $iw_1$ and $iw_2$. By construction
these matrix
coefficients are invariant under interchanging of $w_1$ and $w_2$.
We show duality of this solution $S$ in its spectral variable $\lambda$ and 
its geometric variable.
Consequently it satisfies another two Askey-Wilson second order
difference equations in its spectral variable.

In a series \cite{Ruijs2}, \cite{Ruijs3}, \cite{Ruijs4} of papers, Ruijsenaars
considered a solution $R$ of the same Askey-Wilson difference equations.
These equations arose in his study of relativistic quantum integrable
systems. Ruijsenaars defined the hypergeometric function $R$
 as a Barnes' type integral with integrand expressed in terms of the 
hyperbolic gamma function. Subsequently he established for $R$
duality, $D_4$-symmetry in the parameters, asymptotic behaviour and
the reduction to Askey-Wilson polynomials. We use
these properties to show equality of $R$ to $S$,
which is not apparent from their explicit integral representations.

The structure of this paper is as follows.
In Sections \ref{sec2} and \ref{sec3} we recall some properties of the
hyperbolic gamma function and of Ruijsenaars' hypergeometric function  $R$, 
respectively.
In Section \ref{sec4} we define the modular double of 
$\uq$ and its principal series representation on 
meromorphic functions.
In Section \ref{sec5} we consider the corresponding 
eigenvalue problem of two commuting twisted primitive elements.
Using the
matrix coefficients of the principal series representation 
 we construct a solution $S$ to the Askey-Wilson
difference equations in Section \ref{sec6}, and we 
establish the duality of $S$. In Section \ref{sec7}
we show by a direct calculation that $S$ reduces to the Askey-Wilson
polynomials for certain discrete values of the spectral parameter. Finally, in 
Section \ref{sec8} we show that $S$ equals
 Ruijsenaars' hypergeometric function $R$.

\subsection{Notational conventions}
If $\pm$ appears inside the argument of functions  we mean a product, e.g.
\begin{equation}\label{eqconvention}
f(z\pm a) = f(z+a)f(z-a). 
\end{equation}
Otherwise it means that all sign combinations are possible.

Whenever we use a square root, we always mean the branch which has a cut 
along the negative real line and maps the positive real line to itself.

\subsection{Acknowledgements}
I would like to thank my supervisor Jasper Stokman for many illuminating
conversations and his help in writing this article. 


\section{The hyperbolic gamma function}\label{sec2}
Both Ruijsenaars' and our solution to the 
 Askey-Wilson second order difference equations are expressible in 
terms of the hyperbolic gamma function, which was introduced in
\cite{Ruijs1}. Let us therefore recall some basic properties
of this function, see \cite{Ruijs1} and the appendices of \cite{Ruijs2} 
for more details.

Let us first define for $w_1,w_2 \in \mbb{C}_+ = \{z\in \mbb{C}~|~ \Re(z)>0\}$, 
\begin{equation}\label{eqdefg}
g(w_1,w_2;z) = \int_{0}^{\infty} 
\left( \frac{\sin(2yz)}{2\sinh(w_1 y)\sinh(w_2 y)} - \frac{z}{w_1w_2y}
 \right) \frac{dy}{y}.
\end{equation}
Notice that the integrand has no pole at 0. To ensure convergence
of the integral at infinity however, we must impose the condition
 $|\Im(z)| < \Re(w)$, where $w$ is defined by
\begin{equation*}
w = \frac{w_1+w_2}{2}.
\end{equation*}

The hyperbolic gamma function $G(z)=G(w_1,w_2;z)$ 
for $|\Im(z)| < \Re(w)$ is now defined by
\begin{equation}\label{eqdefG}
G(w_1,w_2;z) = e^{ig(w_1,w_2;z)}.
\end{equation}

The hyperbolic gamma function $G$ owes its name to the fact
that it satisfies the difference equations 
\begin{equation}\label{eqdiffG}
\begin{split}
G(z+iw_1/2) & = 2\cosh(\pi z/w_2) G(z-iw_1/2), \\
G(z+iw_2/2) & = 2\cosh(\pi z/w_1) G(z-iw_2/2).
\end{split}
\end{equation}
In these equations we suppress the $w_1$ and $w_2$ dependence of $G$, 
which we continue to do whenever this does not cause confusion.
These two difference equations allow for an analytic continuation of 
$G$ to a meromorphic function on $\mbb{C}$. The hyperbolic gamma function 
can also be expressed in terms of Barnes' double gamma function,
or Kurokawa's double sine function. Details can be found in 
\cite[Appendix A]{Ruijs2}.

Let us first note a few symmetries of the hyperbolic gamma function, 
which are all obvious from \eqref{eqdefg}:
\begin{align}\label{eqsymG}
   G(w_1,w_2;z) &= G(w_2,w_1;z), \\
   G(w_1,w_2;z) &= G(w_1,w_2;-z)^{-1}, \label{eqsymG2} \\
   G(w_1,w_2;z) &= \overline{G(\bar w_1,\bar w_2,-\bar z)}, \label{eqsymG3} \\
   G(\mu w_1,\mu w_2;\mu z ) & = G(w_1,w_2;z) \qquad (\mu \in \mbb{R}_{>0}). 
 \label{eqsymG4}
\end{align}

The pole and zero locations of $G$  are
easily derived from the difference equations \eqref{eqdiffG}, since
$G$ has no poles or zeros in the strip $z \in \mbb{R} \times
i(-\Re(w),\Re(w))$ in view of \eqref{eqdefG}. The zeros of $G$ are contained
in the set
\begin{equation}\label{eqdeflambdaplus}
\Lambda_+ = iw + iw_1\mbb{Z}_{\geq 0} + iw_2\mbb{Z}_{\geq 0}
\end{equation}
and the poles in $-\Lambda_+$. 
The pole at $z=-iw$ is simple, and its residue equals
\begin{equation}\label{eqresg}
\frac{i}{2\pi}\sqrt{w_1w_2}.
\end{equation}
If $w_1/w_2$ is irrational all other poles are also simple and their residues
 can be calculated from \eqref{eqresg} and
the difference equations \eqref{eqdiffG}, see \cite[Proposition III.3]{Ruijs1}.

For later purposes it is convenient to call an infinite sequence of points in 
$\mbb{C}$ increasing
(respectively decreasing) if it is contained in a set of the form 
$a + \Lambda_+$ (respectively $a-\Lambda_+$) for some $a\in \mbb{C}$. 
In this terminology, $G$ has one increasing zero-sequence and one 
decreasing pole-sequence.

We also need an estimate for $G(z)$ as $\Re(z) \to \infty$ and
$\Im(z)$ stays bounded. In fact we only need it for the quotient of two
hyperbolic gamma functions, which is 
easily derived from the estimate of the hyperbolic gamma itself 
as described in 
\cite[Proposition III.4]{Ruijs1} and \cite[(3.3)]{Ruijs3}.
For $a,b \in \mbb{C}$ and $w_1,w_2 \in (0,\infty)$ the resulting 
estimate reads 
\begin{equation}\label{eqlimitgfrac}
\frac{G(z-a)}{G(z-b)} = 
 \exp\left( \frac{\pi}{2iw_1w_2} (2 z (b-a) + a^2 -b^2 + f(z)) \right),
\end{equation}
where $f(z)$ satisfies for $\Re(z) > \max(w_1,w_2) + \max(\Re(a),\Re(b))$,
\begin{equation}\label{eqpropg}
|f(z)| < C(w_1,w_2,\Im(z),a,b) e^{-\pi\Re(z)/\max(w_1,w_2)},
\end{equation}
with $C$ depending continuously on 
$(0,\infty)^2 \times \mbb{R} \times 
\mbb{C}^2$.

We also use the description of $G$ as a quotient
\begin{equation}\label{eqeandg}
G(z) = \frac{E(z)}{E(-z)},
\end{equation}
where $E$
is an entire function with zeros at $\Lambda_+$
which are all simple if $w_1/w_2$ is irrational.
For a precise definition of $E$, see \cite[Appendix A]{Ruijs2}.

We will occasionally meet functions defined by an integral of the form 
\begin{equation}\label{eqdefm}
M(u,d) = \int_{\mbb{R}} \prod_{j=1}^n 
     \frac{G(w_1,w_2;z-u_j)}{G(w_1,w_2;z-d_j)} dz
\end{equation}
for $w_1,w_2>0$ and for parameters $u_j$ and $d_j$ 
satisfying $|\Im(u_j)|,|\Im(d_j)| <w$ and 
$\Im(\sum_{j=1}^n (u_j-d_j))>0$. These conditions ensure that 
the integral is well defined (the contour meets
no poles and it decreases exponentially at $\pm \infty$).
In \cite[Appendix B]{Ruijs2} it is shown that 
\begin{equation*}
M(u,d) \prod_{j,k=1}^n E(-iw + u_j-d_k)
\end{equation*}
has a unique analytic extension to the set $\{(w_1,w_2,u,d) \in \mbb{C}_+^2 \times 
\mbb{C}^{2n} ~|~ \Im(\sum (u_j-d_j)/w_1w_2)>0\}$. 
Hence $M(u,d)$ is a meromorphic 
function which can only have poles when some $E(-iw +u_j-d_k)$ is zero.


\section{Ruijsenaars' hypergeometric function}\label{sec3}
Ruijsenaars \cite{Ruijs2} introduced a generalization $R$ of the 
 hypergeometric function as a Barnes' type integral.
We recall several properties of $R$ from \cite{Ruijs2} and \cite{Ruijs3}
which we will need to  relate  $R$ to 
the formal matrix coefficients we are going to define in subsequent sections.

We define Ruijsenaars' hypergeometric function
in terms of a parameter set $\gamma_\mu$ $(\mu = 0,1,2,3)$, which
is related to Ruijsenaars' original $c$-parameters by 
\cite[(1.11)]{Ruijs3}. Dual parameters 
$\hat \gamma_\mu$ are defined as
\begin{equation}\label{eqdefhatg}
\left( \begin{array}{c}
\hat \gamma_0 \\ \hat \gamma_1 \\ \hat \gamma_2 \\ \hat \gamma_3
\end{array} \right) = 
\frac{1}{2} \left( \begin{array}{rrrr}
1 & 1& 1 & 1 \\ 1 & 1 & -1 & -1 \\ 1 & -1 & 1 & -1 \\ 1 & -1 & -1 & 1
\end{array} \right) \
\left( \begin{array}{c}
\gamma_0 \\ \gamma_1 \\ \gamma_2 \\ \gamma_3
\end{array} \right) . 
\end{equation}
We denote the set of parameters $(\gamma_0,\gamma_1,\gamma_2,\gamma_3)$
by $\gamma$ and the set of dual parameters by $\hat \gamma$. Note that
taking dual parameters is an involution, $\hat{\hat{\gamma}}=\gamma$.

Ruijsenaars' hypergeometric function $R$ is now defined 
for generic parameters $w_1,w_2 \in \mbb{C}_+, \gamma \in \mbb{C}^4$ by
\begin{multline}\label{eqdefR}
R(\gamma; x,\lambda) = \\
\frac{1}{\sqrt{w_1w_2}} \int_{\mathcal{C}} 
\frac{G(z\pm x+i\gamma_0)G(z\pm \lambda +i\hat \gamma_0)}
       {G(\pm x + i \gamma_0)G(\pm \lambda + i\hat \gamma_0)G(z+iw)}
 \prod_{j=1}^3 \frac{G(i\gamma_0+i\gamma_j + iw)}{G(z+i\gamma_0
        +i\gamma_j + iw)} dz.
\end{multline}
Note that we use the convention \eqref{eqconvention} in this expression.
The integral is taken over a contour $\mathcal{C}$, which is a deformation 
of $\mbb{R}$ seperating the poles of the numerator from the zeros of the 
denominator (equivalently, $\mc{C}$ separates the increasing pole sequences
of the integrand from the downward pole sequences). 
$R$ has an analytic extension to a meromorphic function on 
$(w_1,w_2,\gamma,x,\lambda) \in \mbb{C}_+^2 \times \mbb{C}^6$, with possible 
poles for fixed values of $w_1$, $w_2$, and $\gamma$ at 
\begin{equation}\label{eqpolesR}
 x \in \pm( \Lambda_+ - i\gamma_j ), \qquad
 \lambda \in \pm (\Lambda_+ - i\hat \gamma_j) \qquad (j =0,1,2,3).
\end{equation}
Recall that $\Lambda_+$ is defined
by \eqref{eqdeflambdaplus}. 

We now look at the Askey-Wilson second order difference equations which
$R$ satisfies. The equations are obtained from
\cite[Theorem 3.1]{Ruijs2} by not only replacing the $c$-variables by
$\gamma$, but also multiplying the equations by a constant. 
These descriptions of the Askey-Wilson difference equations
are more convenient for the representation theoretic
approach we consider in the following sections. 

Let us define the function $A$ by
\begin{equation}\label{eqdefar}
\begin{split}
A(w_1,w_2,\gamma;x) &= - \frac{e^{\pi i w_1/w_2+2 \pi i \hat \gamma_0 /w_2}}{
     \sinh\left(\frac{2\pi x}{w_2}\right)
     \sinh\left(\frac{2\pi(x+iw)}{w_2}\right)}
  \prod_{j=0}^3 \cosh\left(\frac{\pi}{w_2}(x+ \frac{iw_1}{2} 
              + i \gamma_j )\right) \\
&=
   \frac{1}{(1-e^{4 \pi x/w_2})(1-e^{4\pi (x+iw)/w_2})} 
\prod_{j=0}^3 (1-e^{2\pi(iw+i\gamma_j +x)/w_2}).
\end{split}
\end{equation}
The Askey-Wilson second order difference operator $\mathcal{L}_\gamma^x$ 
is defined by
\begin{multline}\label{eqdefl}
\mathcal{L}_{\gamma}^x f(x) = 
  A(w_1,w_2,\gamma;x) (f(x+iw_1)-f(x))  \\ +\ A(w_1,w_2,\gamma;-x)  
        (f(x-iw_1)-f(x)) .
\end{multline}
Here the superscript $x$ is added to emphasize that the operator acts on the
$x$-variable (in a moment we will also consider the operator $\mc{L}$ acting
on the spectral variable $\lambda$).
We write $\widetilde{\mc{L}}_{\gamma}^{x}$ for the Askey-Wilson operator
\eqref{eqdefl} with $w_1$ and $w_2$ interchanged. 

Ruijsenaars' hypergeometric function $R$ is an eigenfunction of four  
Askey-Wilson second order difference operators with eigenvalues 
expressible in terms of
\begin{multline}\label{eqdefv}
v(w_1,w_2,\gamma;\lambda) =  \\
-2e^{\pi i w_1/w_2+ 2\pi i \hat \gamma_0 /w_2} (
\cosh(2\pi\lambda/w_2) + \cosh(\pi i w_1/w_2 + 2 \pi i \hat \gamma_0/w_2)).
\end{multline}
Specifically, $R$ satisfies the
difference equations
\begin{equation}\label{eqrdiff}
\begin{split}
\mc{L}_\gamma^x R(\gamma;x,\lambda) &= 
   v(w_1,w_2,\gamma;\lambda) R(\gamma;x,\lambda), \\
\tc{L}_\gamma^x R(\gamma;x,\lambda) & = 
   v(w_2,w_1,\gamma;\lambda) R(\gamma;x,\lambda), \\
\mc{L}_{\hg}^\lambda R(\gamma;x,\lambda) &= 
   v(w_1,w_2,\hg;x) R(\gamma;x,\lambda), \\
\tc{L}_{\hg}^\lambda R(\gamma;x,\lambda) &=
   v(w_2,w_1,\hg;x) R(\gamma;x,\lambda).
\end{split}
\end{equation}
Actually the last three of these equations follow from the first
by various symmetries of $R$. The second 
difference equation can be derived from the first
(and the fourth from the third) by using the fact that 
$R$ is invariant under the exchange of $w_1$ and $w_2$,
\begin{equation*}
R(w_1,w_2,\gamma;x,\lambda) = R(w_2,w_1,\gamma;x,\lambda).
\end{equation*}
This symmetry can be directly seen from the definition \eqref{eqdefR} of $R$ 
 and the corresponding symmetry \eqref{eqsymG} of the hyperbolic gamma
function.
The third difference equation can be obtained from the first by using the 
duality of $R$ under the exchange of $x$ and $\lambda$, 
\begin{equation}\label{eqrdual}
R(\gamma;x,\lambda) = R(\hg;\lambda,x).
\end{equation}
This duality is also a direct consequence of the definition of $R$ 
using the fact that $\gamma_0 + \gamma_j = \hg_0 + \hg_j$ for $j=1,2,3$.

There are more symmetries of $R$ 
directly visible from the definition. Since the hyperbolic gamma
function is scale invariant it follows that $R$ is scale invariant as well,
\begin{equation*}
R(\nu w_1,\nu w_2, \nu \gamma;\nu x, \nu \lambda) = 
  R(w_1,w_2,\gamma;x,\lambda)
\end{equation*}
for $\nu \in (0,\infty)$, where $\nu \gamma$ denotes the scaled parameter set
$(\nu \gamma_0,\nu \gamma_1,\nu \gamma_2, \nu \gamma_3)$. 
Furthermore it is
immediately clear that $R$ is symmetric under 
permutations of $\gamma_1$, $\gamma_2$, and $\gamma_3$. This symmetry can
be extended to a $D_4$-symmetry in the four parameters $\gamma$ 
(where the Weyl
group of type $D_4$ acts on the parameters by permutations and an even number
of sign flips). To formulate this result we need the $c$-function
\begin{equation*}
c(\gamma;y) = \frac{1}{G(2y+iw)} \prod_{j=0}^3 G(y-i\gamma_j)
\end{equation*}
and the normalization constant
\begin{equation}\label{eqdefn}
N(\gamma) = \prod_{j=1}^3 G(i\gamma_0+ i\gamma_j+iw).
\end{equation}
The $D_4$-symmetry \cite[Theorem 1.1]{Ruijs3} of $R$ then reads
\begin{equation}\label{eqd4sym}
\frac{R(\gamma;x,\lambda)}{c(\gamma;x)c(\hat \gamma;\lambda) N(\gamma)}
 = \frac{R(w(\gamma);x,\lambda)}
   {c(w(\gamma);x)c(\widehat{w(\gamma)};\lambda)N(w(\gamma))}
\end{equation}
 for all elements $w$ of the Weyl group of type $D_4$.
Notice that both the $c$-function and $N$ are invariant under the action 
of the $S_3$-subgroup which permutes $\gamma_1$, $\gamma_2$, and $\gamma_3$.

Finally we recall the limit behaviour of $R$, cf.\ \cite[Theorem 1.2]{Ruijs3}. 
Set $\alpha = 2\pi/w_1w_2$. For $w_1,w_2 \in \mbb{R}_{>0}$, 
$\gamma \in \mbb{R}^4$, and $w_1\neq w_2$ there exists an open
neighbourhood $U\subset \mbb{C}$ of $\mbb{R}$, such that
the asymptotics of $R$ for fixed $\lambda \in U$ are given by
\begin{equation}\label{eqRasymp}
R(\gamma;x,\lambda) =   
\mc{O} (e^{\alpha(|\Im(\lambda)| - \hg_0 - w)|\Re(x)|})
\end{equation}           
for $\Re(x) \to \pm \infty$, uniformly for $\Im(x)$ in compacta. In fact,
Ruijsenaars gives a precise expression for the leading term of $R$ as
 $\Re(x)\to \pm \infty$ when $\lambda \in \mbb{R}$.
These results easily extend to $\lambda$ in some
open neighbourhood $U$ of $\mbb{R}$.                    


\section{The modular double of $\uq(\mathfrak{sl}_2(\mbb{C}))$}\label{sec4}
In this section we consider a slightly extended version of Faddeev's 
\cite{Faddeev} modular double
of $\mathcal{U}_q(\mathfrak{sl}_2(\mathbb{C}))$ and define an algebraic 
version of its principal series 
representation on the space $\mc{M}$ of meromorphic functions on $\mbb{C}$. 
We define an inner product on some suitable subspace of $\mc{M}$, which is 
compatible to the $*$-structure on $\uq(\mathfrak{sl}_2(\mbb{C}))$ 
associated to the real form $\mathfrak{sl}_2(\mbb{R})$ of 
$\mathfrak{sl}_2(\mbb{C})$, cf.\ \cite{Masuda}.

Throughout Sections \ref{sec4}, \ref{sec5} and \ref{sec6} 
we assume that  $w_1$ and $w_2$ are positive
real numbers such that neither
$w_1/w_2$ nor $w_2/w_1$ is an integer, unless specifically 
stated otherwise. We define
\begin{equation*}
q = \exp(\pi i w_1/w_2), \qquad \tilde q = \exp(\pi i w_2/w_1),
\end{equation*}
which both lie on the unit circle (but they are not $\pm 1$).
For complex numbers $\beta$ we define
\[
q^{\beta} = e^{\beta \pi i w_1/w_2}, \qquad
\tilde q ^{\beta} = e^{\beta \pi iw_2/w_1}.
\]

\begin{defi}
The quantized universal enveloping algebra 
 $\mc{U}_q = \mc{U}_q(\mathfrak{sl}_2(\mathbb{C}))$ of 
$\mathfrak{sl}_2(\mbb{C})$ is the unital associative algebra 
over $\mbb{C}$ generated by $K^{\pm 1}$, $E$, and $F$, subject to the 
relations
\begin{equation*}
\begin{split}
& KK^{-1}=K^{-1}K=  1, \\
& KE = q^2 EK, \\ & KF = q^{-2}FK, \\
& EF-FE = \frac{K-K^{-1}}{q-q^{-1}}.
\end{split}
\end{equation*}
\end{defi}
If $w_1/w_2$ is irrational, then 
the center of $\uq$ is generated by the quantum Casimir element $\Omega$,
defined as
\begin{equation*}
\Omega = qK + q^{-1}K^{-1} + (q-q^{-1})^2 FE.
\end{equation*}

By simply replacing $q$ by $\tilde q$ (or interchanging $w_1$ and $w_2$) we
obtain the quantum universal enveloping algebra
$\mc{U}_{\tilde q}$. The generators of $\uqt$
are denoted by $\tilde K^{\pm 1}$, $\tilde E$, and $\tilde F$.
The following concept of modular double was introduced by Faddeev
\cite{Faddeev}.
\begin{defi} The modular double $\mc{Q}$ is $\uq \otimes \uqt$ 
endowed with its standard tensor product algebra structure.
\end{defi}
For elements $X \in \uq$ (respectively $\tilde X \in \uqt$) 
we also write $X$ (respectively $\tilde X$) for its image under
the natural embedding of $\uq$ (respectively $\uqt$) in $\mc{Q}$. 
In particular, $X \tilde X=\tilde X X$ in $\mc{Q}$ for  elements 
$X\in \uq$ and $\tilde X \in \uqt$.

We now define an extension of the modular double by formally adjoining 
complex powers of $K$ and $\tilde K$ to $\mc{Q}$.
Let $\mc{A} = \bigoplus_{x\in \mbb{C}}\mbb{C}\hat{x}$ 
be the group algebra
of the additive group $\hat{\mbb{C}} = (\mbb{C},\oplus)$, 
where $\oplus$ is the
translated addition $\hat x \oplus \hat y = 
\widehat{x\!\!+\!\!y\!\!+\!\!iw}$ (this translation in addition will
make formules simpler later on). The unit of
$\hat{\mbb{C}}$ is $\widehat{-iw}$.

\begin{lemma}\label{lemaction}
There exists a unique left $\mc{A}$-action by
algebra automorphisms on the modular double $\mc{Q}$ satisfying  
\begin{align*} 
&\hat x \cdot K^{\pm 1} = K^{\pm 1},  &
&\hat x \cdot \tilde K^{\pm 1} = \tilde K^{\pm 1}, \\
&\hat x \cdot E = -qe^{2\pi x/w_2} E, &
&\hat x \cdot \tilde E = -\qt e^{2\pi x/w_1} \tilde E, \\
&\hat x \cdot F = -qe^{-2\pi x/w_2} F,  &
&\hat x \cdot \tilde F = -\qt e^{-2\pi x/w_1} \tilde F.
\end{align*}
\end{lemma}
\begin{proof}
Observe that e.g. the action of $\hat x$ on $E$ can be rewritten as
\[
\hat x \cdot E = e^{2\pi (x+iw)/w_2} E.
\]
The lemma now follows by direct calculations.
\end{proof}

\begin{defi} The extended modular double $\mathcal{D} 
= \mc{Q} \rtimes \mc{A}$ 
is the crossed product of the modular double $\mc{Q}$ and the
algebra $\mc{A}$ under its action on $\mc{Q}$ as defined
in Lemma \ref{lemaction}. 
\end{defi}
Hence $\mc{D}$ is the vector space 
$\mc{Q} \otimes \mc{A}$ endowed  with the unique algebra structure 
such that the natural embeddings of 
$\mc{Q}$ and $\mc{A}$ in $\mc{D}$ are algebra morphisms and such that
\[
\hat x Q = (\hat x \cdot Q) \hat x, \quad \forall x \in \mbb{C},\; \forall 
Q\in \mc{Q},
\]
where we identified $\hat x$ (respectively $Q$) with their
images under the natural embeddings of $\mc{A}$ (respectively
$\mc{Q}$) in $\mc{D}$.

Now we define representations $\phl$ of the extended modular double $\mc{Q}$ 
on the space $\mc{M}$ of meromorphic functions on $\mbb{C}$ depending on 
a complex representation label $\lambda$, cf. \cite{Kharchev}. 
These representations may be viewed as algebraic versions of the 
principal series 
representations of $\uq(\mathfrak{sl}_2(\mbb{R}))$. We define these
representations in terms of the operators $T_y$ and $S_y$ on $\mathcal{M}$,
which act by
\[
T_y f(z) = f(z+y), \qquad S_y f(z) = e^{2\pi i z/y}f(z) \qquad
 (y\in\mbb{C}).
\]
\begin{lemma} \label{lemrep}
For $\lambda \in \mbb{C}$ the assignments
\begin{align*}
\phl(K) &= T_{iw_1},\quad  \phl(\tilde K) = T_{iw_2}, \quad \phl(\obx) = T_{x+iw}, \\
\phl(E) &= \frac{q^{1/2}}{q-q^{-1}} S_{iw_2} 
    \left( q^{-1/2}e^{\pi\lambda/w_2} + q^{1/2}e^{-\pi\lambda/w_2} T_{iw_1}
     \right), \\
\phl(F) &= -\frac{q^{1/2}}{q-q^{-1}} S_{-iw_2} 
    \left( q^{-1/2}e^{\pi\lambda/w_2} + q^{1/2}e^{-\pi\lambda/w_2} T_{-iw_1}
  \right), \\
\phl(\tilde E) &= \frac{\qt^{1/2}}{\qt-\qt^{-1}} S_{iw_1} 
    \left( \qt^{-1/2}e^{\pi\lambda/w_1} + 
     \qt^{1/2}e^{-\pi\lambda/w_1} T_{iw_2} \right), \\
\phl(\tilde F) &= -\frac{\qt^{1/2}}{\qt-\qt^{-1}} S_{-iw_1} 
    \left( \qt^{-1/2}e^{\pi\lambda/w_1} + 
    \qt^{1/2}e^{-\pi\lambda/w_1} T_{-iw_2} \right) ,
\end{align*}
uniquely define a representation $\phl$ of $\mc{D}$ on $\mc{M}$.
\end{lemma}
Observe that the action of the generators of $\uqt$ are obtained from the 
action of the generators of $\uq$ by interchanging $w_1$ and $w_2$.
\begin{proof}
The defining relations of $\mc{D}$ are easily checked  using
$T_xT_y = T_{x+y} = T_yT_x$, $S_xS_{-x} = 1$, and the equation
\begin{equation*}
T_x S_y = e^{2\pi i x/y}  S_y T_x.
\end{equation*}
\end{proof}

\begin{remark}
Denote $v=(w_1-w_2)/2$, then $\phl(\widehat{iv})=\phl(K)$ and
$\phl(\widehat{-iv}) = \phl(\tilde K)$.
The extension of the modular double $\mc{Q}$ by $\mc{A}$ and the extension
of the representation $\left. \phl \right|_{\mc{Q}}$ to $\phl$ thus have the
effect of introducing
non-integral powers of $T_{iw_1}$ and $T_{iw_2}$ in the image of $\phl$.
The introduction of this extension is not an essential part of the analysis 
later on and
is only included for simplification. Using only integral powers of 
$K$ and $\tilde K$ we can simulate the action of $\hat x$ for $x$ in some
dense subset of $\mbb{R}$, cf.\ \cite[Proposition 1.6]{Kharchev}.
\end{remark}

A simple calculation shows that $\phl(\Omega)$ acts as a scalar,
\begin{equation}\label{eqrepomega}
\phl(\Omega) f = -2\cosh(2\pi\lambda/w_2) f, \qquad f\in \mc{M}.
\end{equation}
Since $\phl$ is an algebraic
version of the  principal series representation with 
representation label $\lambda \in \mbb{C}$, this is as expected. 

\begin{defi} 
We say that $f\in \mc{M}$ 
has exponential growth with growth rate $\epsilon \in \mbb{R}$
if there exists a
compact set $K_f \in \mbb{R}$  such that all poles of $f$ are contained
in $K_f \times i\mbb{R} = \{ x+iy ~|~ x \in K_f, y \in \mbb{R}\}$ and
if $|f(x+iy)| = \mc{O}(\exp(\epsilon|x|))$ for $x \to \pm \infty$, 
uniformly for $y$ in compacta of $\mbb{R}$.
\end{defi}

On the space of meromorphic functions which have negative
exponential growth and which have no poles
on $\mbb{R}$, we define a sesquilinear form by
\begin{equation}\label{eqdefip}
\langle f,g \rangle = \int_{\mbb{R}} f(z)\overline{g(z)}dz.
\end{equation}
Observe that this expression is already well defined under the milder 
asymptotic condition that the sum
of the two exponential growths of $f$ and $g$ is negative. 
Note furthermore that \eqref{eqdefip} can be rewritten as
\begin{equation}\label{eqdefip2}
\langle f,g \rangle = \int_{\mbb{R}} f(z)\bar g(z) dz,
\end{equation}
where $\bar g(z):= \overline{g(\bar z)}$ now is analytic
at  $z \in \mbb{R}$.

Following \cite{Kharchev}
we define an antilinear anti-algebra involution $*$ on the
extended modular double $\mc{D}$ by
\begin{equation}\label{eqstarstructure}
K^* = K,\; E^* = -E,\; F^*=-F,\; \tilde K^* = \tilde K,\; 
\tilde E^* = -\tilde E,\; \tilde F^* = - \tilde F,\; 
\obx^* = \widehat{-\bar x}.
\end{equation}
If we restrict this involution to $\uq$ (respectively $\uqt$)
we obtain the $*$-structure on $\uq$ (respectively $\uqt$) 
corresponding to the noncompact real form 
$\mathfrak{sl}_2(\mbb{R})$ of $\mathfrak{sl}_2(\mbb{C})$, cf.\ \cite{Masuda}.

The following lemma relates the sesquilinear form \eqref{eqdefip} to the 
$*$-structure \eqref{eqstarstructure} on $\mc{D}$.
\begin{lemma}\label{lemadjoint}
Let $\lambda \in \mbb{C}$ and $f,g \in \mc{M}$. If 
the poles of $f$ and $g$ are outside the strip $\mbb{R} \times 
i[-w_1,w_1]$ and if the sum of the exponential growth rates of $f$ and $g$ is
smaller than $ -2\pi/w_2$, then 
\[
\langle \phl(X)f,g \rangle = \langle f,\phlb(X^*)g \rangle
\]
for 
$X \in \mc{U}_{q,1} := {\mbox{span}}_{\mbb{C}}\{1,E,F,K,K^{-1},FK,EK^{-1}\}$.
\end{lemma}
\begin{proof}
In view of \eqref{eqdefip2}
the proof follows by  a change of variables and some contour shifting using
Cauchy's theorem.
\end{proof}

A similar lemma holds for the dual algebra $\uqt$.


\section{Twisted primitive elements and matrix coefficients}\label{sec5}
Koornwinder \cite{Koorn} introduced twisted primitive elements to
obtain the Askey-Wilson polynomials as matrix coefficients 
 of finite dimensional 
$\mathcal{U}_q(\mathfrak{sl}_2(\mbb{C}))$-re\-pre\-sen\-ta\-tions.
We recall the definition of twisted primitive elements 
and  show that they act as first order
difference operators under the representations $\phl$. 
We construct eigenvectors to these operators in
terms of the hyperbolic gamma function and we consider the corresponding
formal matrix coefficients $\psi$ of $\phl$. In subsequent sections we 
relate $\psi$ to Ruijsenaars' hypergeometric function.

Let $\rho\in \mbb{C}$ and set
\begin{equation*}
\nu_\rho = q^{2\rho/w_1} + q^{-2\rho/w_1} = 2\cos(2\pi \rho/w_2).
\end{equation*}
The twisted primitive element 
$Y_\rho \in \uq \subset \mathcal{D}$ is defined as 
\begin{equation}\label{eqdefyrho}
Y_\rho = iq^{-1/2}E + i q^{-1/2} FK - \frac{\nu_\rho}{q-q^{-1}} (K-1).
\end{equation}

Analogously we define the twisted primitive element $\tilde Y_\rho \in \uqt$ 
by interchanging $w_1$ and $w_2$, viz.
\begin{equation*}
\tilde Y_\rho = i \qt^{-1/2}\tilde E + i \qt^{-1/2} \tilde F \tilde K
    - \frac{\tilde \nu_\rho}{\qt-\qt^{-1}} (\tilde K-1),
\end{equation*}
where $\tilde \nu_\rho = 2\cos(2\pi  \rho/w_1)$. 

Denoting
\begin{equation*}
\mu_\tau(\rho) = \frac{\nu_\rho- \nu_\tau}{q-q^{-1}}, \quad \tilde 
\mu_{\tau}(\rho)  = \frac{\tilde \nu_\rho- \tilde \nu_\tau}{\qt-\qt^{-1}},
\end{equation*}
we now have the following lemma.
\begin{lemma}\label{lemhtr}
The function 
\[
H_{\tau,\rho}^\lambda(z) = 
\frac{G(z+\lambda/2-3iw/2 \pm i\tau)}
     {G(z-\lambda/2-iw/2 \pm i\rho)}
\]
satisfies
\begin{equation}\label{eqtempq}
\begin{split}
\phl(Y_{\rho}) H_{\tau,\rho}^\lambda & = 
     \mu_{\tau}(\rho) H_{\tau,\rho}^\lambda, \\
\phl(\tilde Y_{\rho}) H_{\tau,\rho}^\lambda &=
     \tilde \mu_{\tau}(\rho) H_{\tau,\rho}^\lambda.
\end{split}
\end{equation}
\end{lemma}
\begin{proof}
Since $\htr$ is invariant under the exchange of $w_1$ and
$w_2$, it is sufficient to prove only the first eigenvalue equation.
A calculation shows that $\phl(Y_\rho)f=\mu_\tau(\rho)f$ is
equivalent to the first order difference equation
\begin{equation}\label{eqfode}
f(z+iw_1/2) = 
\frac{\cosh(\frac{\pi}{w_2}(z+\lambda/2-3iw/2\pm i\tau))}
 {\cosh(\frac{\pi}{w_2}(z-\lambda/2 -iw/2 \pm i\rho))} f(z-iw_1/2).
\end{equation}
(The exact calculation can be found in Appendix \ref{app1}.)
Using the difference equation \eqref{eqdiffG} for the hyperbolic gamma 
function it immediately follows that $\htr$ satisfies the difference equation 
\eqref{eqfode}.
\end{proof} 

\begin{remark}\label{remark2}
For any one of the two
equations \eqref{eqtempq} 
there are infinitely many solutions (we can e.g. multiply a 
solution to the first equation by any $iw_1$-periodic function). 
The crucial step in finding common solutions to both difference 
equations is to rewrite the first
difference equation in the specific form \eqref{eqfode}.
Indeed, the resulting solution $\htr$ in terms of hyperbolic gamma functions 
is invariant under interchanging $w_1$ and $w_2$, hence it automatically 
satisfies the second difference equation.
This is the main difference between our analysis and the one in \cite{Stokman}.
\end{remark}

Now let us consider the adjoint $Y_\sigma^*$, which is
\[
Y_\sigma^* = iq^{1/2}E + iq^{-3/2}FK + \frac{\nu_{\bar \sigma}}{q-q^{-1}}
 (K-1).
\]
Since $\mu_{\upsilon}(\sigma)^* = -\mu_{\bar \upsilon}(\bar \sigma)$, we 
are interested in solutions to the equation 
$\pi_{\bar \lambda}(Y_\sigma^*) f= - \mu_{\bar \upsilon}(\bar \sigma)f$ and 
the corresponding equation $\pi_{\bar \lambda}(\tilde Y_\sigma^*) f = 
- \tilde \mu_{\bar \upsilon}(\bar \sigma) f$
for the second component of the modular double.
\begin{lemma}
The function 
\[
\fus(z) = \frac{G(z + \bar \lambda -iw/2 \pm i \bar \upsilon)}{
           G(z- \bar \lambda + iw/2 \pm i \bar \sigma)}
\]
satisfies
\begin{align*}
\pi_{\bar \lambda}(Y_{\sigma}^*) \fus & = -\mu_{\bar \upsilon}(\bar \sigma) 
 \fus,\\
 \pi_{\bar \lambda}(\tilde Y_{\sigma}^*) \fus & = 
          -\tilde \mu_{\bar \upsilon}(\bar \sigma) \fus.
\end{align*}
\end{lemma}
\begin{proof}
The proof is similar to the proof of the previous lemma.
\end{proof}

We will need a few results on the analytic properties
of the two functions $H^{\lambda}_{\tau,\rho}$ and $F^{\lambda}_{\upsilon,
\sigma}$.

\begin{lemma}
The possible pole locations of $\htr$ and $\fus$ are  at
\begin{equation*}
  -\lambda/2 \pm i \tau + iw - \Lambda_+,
          \qquad \lambda/2 \pm i \rho + iw + \Lambda_+
\end{equation*}
and 
\begin{equation*}
   -\bar \lambda/2 \pm i \bar \upsilon - \Lambda_+, \qquad
  \bar \lambda/2 \pm i \bar \sigma + \Lambda_+
\end{equation*}
respectively. Furthermore, $\htr$ and $\fus$ have
 exponential growth with growth rates 
$\pi ( 2\Im(\lambda)-2w)/w_1w_2$ and 
$\pi (- 2 \Im(\lambda)-2w)/w_1w_2$, respectively.
\end{lemma}
\begin{proof}
The proof follows directly from the zero/pole locations and
asymptotics of the 
hyperbolic gamma function (see Section \ref{sec2}).
\end{proof}

Define
\begin{equation}\label{eqdefxi}
\xi = \max(|\Re(\rho)|,|\Re(\sigma)|,|\Re(\tau)|,|\Re(\upsilon)|)
\end{equation}
and
\begin{equation}\label{eqdefzeta}
\zeta = w/2 - \xi - |\Im(\lambda/2)|.
\end{equation}
We assume that the parameters $\rho,\sigma,\upsilon,\tau$ and
the variable $\lambda$ are such that $\zeta>0$.
For $|\Im(x)| <\zeta$ define
\begin{equation}\label{eqdefpsi}
\psi(\rho,\sigma,\tau,\upsilon;x,\lambda) = \langle \phl(\hat x)\htr,\fus \rangle,
\end{equation}
which is well defined since the exponential growth of the integrand 
 equals
$-2\pi (w_1+w_2) /w_1w_2 <0$ and the pole sequences  of $\phl(\hat x)\htr $ and
$\fus$ all stay away from the real line due to the condition $|\Im(x)|<\zeta$.
 Note that the increasing pole sequences of 
$\phl(\hat x)\htr$ and $\fus$ are
all located above the real line and the decreasing pole sequences are all
located below the real line due to the shifted addition
in $\hat{\mathbb{C}}$. 
Observe furthermore that the matrix coefficient $\psi$ is 
invariant under the exchange of $w_1$ and $w_2$, cf.\ Remark \ref{remark2}.
The function $\psi$ will
be related to Ruijsenaars' hypergeometric function $R$ in Section \ref{sec8}.

Using \eqref{eqsymG2}, \eqref{eqsymG3}, and \eqref{eqdefip2} 
we can write $\psi$ as
\begin{equation}\label{eqintpsi}
\psi(\rho,\sigma,\tau,\upsilon;x,\lambda) = 
\int_{\mbb{R}} \frac{G(z+x+\lambda/2-iw/2 \pm i\tau)
          G(z- \lambda/2 -iw/2 \pm i \sigma)}
     {G(z+x-\lambda/2 + iw/2 \pm i \rho)
          G(z + \lambda/2 + iw/2 \pm i \upsilon)} dz,
\end{equation}
which is of the form \eqref{eqdefm}. It follows from the discussion at the
end of Section \ref{sec2} that
\begin{equation}\label{eqdefPsi}
\begin{split}
\Psi (\gamma;x,\lambda) & =  E(x \pm i \gamma_0) E(x \pm i \gamma_1) E(-x \pm i \gamma_2)
  E(-x \pm i \gamma_3)\\
 & \qquad \times   E(\lambda \pm i \hat \gamma_0)
   E(\lambda \pm i \hat \gamma_1)
  E(-\lambda \pm i \hat \gamma_2) 
   E(-\lambda \pm i \hat \gamma_3) \psi(\gamma;x,\lambda)
\end{split}
\end{equation}
has an entire extension to 
\begin{equation}\label{eqdefO}
\mc{O} =  \{ (w_1,w_2,\rho,\sigma,\tau,\upsilon,x,\lambda) 
\in \mbb{C}_+^2 \times \mbb{C}^6 \}.
\end{equation}
Hence $\psi$ is  meromorphic on the same domain $\mc{O}$.


\section{The Askey-Wilson difference equations}\label{sec6}

We show that the formal matrix coefficient
$\psi$ (see Section \ref{sec5}) satisfies a second order difference
equation with step size $iw_1$
using a radial part calculation of the Casimir $\Omega$ with 
respect to twisted primitive
elements. As a consequence a renormalization $S$ \eqref{eqdefs} of
 $\psi$ satisfies an Askey-Wilson second order difference equation. 
Since $S$, like $\psi$, is invariant under exchanging $w_1$ and $w_2$, we
obtain a second difference equation with step size $iw_2$. 
We furthermore show that $S$ satisfies a duality in the
geometric and spectral variables, and we derive various obvious
symmetries of $S$.

Let us start by establishing a correspondence between the set of parameters 
$\rho$, $\sigma$, $\tau$, and $\upsilon$ and Ruijsenaars' parameter set 
$\gamma$ by
\begin{equation}\label{eqmy2ruijs}
\gamma_0 = -\rho + \sigma, \quad
\gamma_1 = \rho + \sigma, \quad
\gamma_2 = - \tau - \upsilon, \quad
\gamma_3 = \tau- \upsilon.
\end{equation}
Observe that  $\hat \gamma_0$ (see 
\eqref{eqdefhatg}) becomes
\[
\hat \gamma_0 = \frac{1}{2}(\gamma_0 + \gamma_1 + \gamma_2 + \gamma_3) = 
          \sigma-\upsilon.
\]
We will also use the abbreviation $\gamma$ for the parameters 
$(\rho,\sigma,\tau,\upsilon)$. In particular we write 
$\psi(\gamma;x,\lambda)$
for \eqref{eqdefpsi}. Later we show that \eqref{eqmy2ruijs} is
the parameter correspondence which relates $\psi$ to Ruijsenaars'
hypergeometric function. 

Now we perform a radial part calculation of the Casimir element 
$\Omega$ with respect to 
the twisted primitive elements (see \eqref{eqdefyrho}).
The result is stated in terms of the function $A$, see \eqref{eqdefar}. 
\begin{lemma}\label{lemradialpart}
We have 
\begin{equation*}
\obx \Omega = \obx \Omega(x) \mod
     \obx \mc{U}_{q,2}(Y_\rho-\mu_{\tau}(\rho)) + 
     (Y_\sigma-\mu_{\upsilon}(\sigma)) \obx \mc{U}_{q,2},
\end{equation*}
where 
$\mc{U}_{q,2} :={\mbox{span}}_{\mbb{C}}\{1,K^{-1},F\}$
and $\Omega(x)$ is defined as the element
\begin{equation*}
\Omega(x)  =  B(x) K + C(x) + D(x)K^{-1}
\end{equation*}
with coefficients
\begin{align*}
B(x) &  = q^{-1}A(\gamma_0, -\gamma_0, \gamma_1, -\gamma_1;x),  \\
C(x) &  =  q^{-1 + 2\hat \gamma_0/w_1}  
        \left[ -A(\gamma;x)-        A(\gamma;-x)   + 1  + 
             q^{2-4\hat \gamma_0/w_1} \right],  \\
D(x) &  = q^{-1}A(\gamma_2,-\gamma_2,\gamma_3,-\gamma_3;-x). 
\end{align*}
\end{lemma}
\begin{proof}
The proof involves a radial part calculation similar to the
one performed in \cite[Proposition 3.3]{Stokman}. In fact we can use 
the calculation in \cite{Stokman}
 using an embedding $\phi$ of the extended quantum universal enveloping
algebra
$\uq \rtimes \hat{\mbb{C}}$ into the one in \cite{Stokman}, given by
\begin{align*}
\phi(K) & = K^2,   & \phi(\hat x) &= 
    \widehat{\left(\frac{2x\!\!+\!\!iw_2}{iw_1}\right)} K, \\
\phi(E) & = -iKX^+,  &  \phi(F) & = iX^-K^{-1}.
\end{align*}
A direct calculation gives $\phi(\Omega) = (q-q^{-1})^2\Omega +2$ and
$\phi(Y_{\rho}) = Y_{2\rho/w_1}$
(on the right hand side we use the $\Omega$ and $Y$ from \cite{Stokman}, 
which have a slightly different definition).
Note that in \cite{Stokman} the radial part is calculated modulo a 
larger vector space. However, it is easily verified that 
the present smaller space suffices for the proof.
\end{proof}

Using this radial part calculation we can prove that $\psi$ \eqref{eqdefpsi}
satisfies a gauge transformed Askey-Wilson second order difference equation.
\begin{lemma}\label{lemsecdifeq}
The function $\psi(x) = \psi(\gamma;x,\lambda)$ 
satisfies the difference equation
\begin{equation}\label{qwerty}
-2\cosh(2\pi\lambda/w_2)\psi(x) = 
   B(x)\psi(x+iw_1) + C(x)\psi(x) + 
  D(x)\psi(x-iw_1),
\end{equation}
and a similar equation with $w_1$ and $w_2$ interchanged. These equations
hold as identities between meromorphic functions on the domain
$\mc{O}$ \textup{(}see \eqref{eqdefO}\textup{)}.
\end{lemma}
\begin{proof}
Observe that by the symmetry of $\psi$ in $w_1$ and $w_2$ we only
have to prove the difference equation \eqref{qwerty}.

We first prove the lemma under restricted parameter conditions, 
which allow us to use the expression \eqref{eqdefpsi} 
of $\psi$ as a matrix coefficient of 
the $\mc{D}$-representation $\phl$. 
Using analytic continuation we can subsequently
remove these parameter constraints, cf.\ the discussion at the end of 
Section \ref{sec5}.

Let us assume that $w_1,w_2>0$ and that
\begin{equation}\label{eqweenie}
w_2>7w_1 + 4\xi + 2|\Im(\lambda)| + 4|\Im(x)| 
\end{equation}
holds. Then $|\Im(x)| < \zeta$, so $\psi$ is defined
by \eqref{eqdefpsi} (recall that  $\xi$ and $\zeta$ are defined by
\eqref{eqdefxi} and \eqref{eqdefzeta}, respectively).
By \eqref{eqrepomega},
\begin{equation}\label{eqstep1}
-2\cosh(2\pi\lambda/w_2) \psi(x) = \langle \phl(\obx \Omega)\htr,\fus \rangle
\end{equation}
holds. By Lemma \ref{lemradialpart} there exist $X,Z\in \mc{U}_{q,2}$ such that
\begin{equation}\label{eqsht1}
\obx \Omega = \obx \Omega(x) + \obx X (Y_\rho - \mu_{\tau}(\rho))
                 +(Y_\sigma - \mu_{\upsilon}(\sigma))\obx Z.
\end{equation}
Since $\phl(Y_{\rho}- \mu_{\tau}(\rho))\htr = 0$, we have
\begin{equation}\label{eqsht2}
 \langle \phl(\obx X(Y_\rho - \mu_{\tau}(\rho))) \htr,\fus \rangle = 0.
\end{equation}
The exponential growth of $\phl(\hat x Z)\htr$ is at most the
exponential growth of $\htr$ plus $2\pi/w_2$ (due to the possible
occurrence of an $S_{iw_2}$ factor in $\pi_{\lambda}(Z)$). 
The sum of the exponential 
growths of $\phl(\hat x Z)\htr$ and $\fus$ is at most $-2\pi/w_1$, hence
strictly smaller than $-2\pi/w_2$, since the restrictions on the parameters 
imply that $w_2>w_1$. Moreover the condition \eqref{eqweenie} implies 
 that neither $\phl(\hat x Z) \htr$ nor $\fus$ has any
poles in the strip $\mbb{R} \times i[-w_1,w_1]$. Using Lemma
\ref{lemadjoint} and the
fact that $Y_\sigma \in \mc{U}_{q,1}$, we thus obtain
\begin{equation}\label{eqsht3}
\langle \phl((Y_{\sigma}-\mu_{\upsilon}(\sigma))\obx Z)
 H_{\tau,\rho},F_{\upsilon,\sigma} \rangle =
 \langle \phl
(\obx Z)H_{\tau,\rho},\phlb(Y_{\sigma}^* + \mu_{\bar \upsilon}(\bar \sigma))
    F_{\upsilon,\sigma} \rangle =0 .
\end{equation}
Combining \eqref{eqstep1},
\eqref{eqsht1}, \eqref{eqsht2}, and \eqref{eqsht3} now yields
\begin{equation}\label{eqbijna}
-2\cosh(2\pi\lambda/w_2) \psi(x) = \langle\phl(\hat x \Omega)\htr,\fus \rangle = 
 \langle \phl(\hat x \Omega(x)) \htr,\fus \rangle .
\end{equation}
Furthermore, by Lemma \ref{lemradialpart}  (remember that 
$\widehat{x\!+\!iw_1}$ and $\obx K$ act in the same way under $\phl$) we have
\begin{equation}\label{eqsht4}
 \langle \phl(\obx \Omega(x))H_{\tau,\rho}, F_{\upsilon,\sigma}\rangle = 
B(x)\psi(x+iw_1) +C(x) \psi(x) +D(x)\psi(x-iw_1). 
\end{equation}
The lemma for the restricted parameter conditions follows now directly
from \eqref{eqbijna} and \eqref{eqsht4}. 
\end{proof}

Using the function 
\[
\Delta(\gamma;x) = \frac{G(x+i\gamma_2)G(x+i\gamma_3)}
                        {G(x-i\gamma_0)G(x-i\gamma_1)},
\]
we can define a renormalization $S$ of $\psi$ as
\begin{equation}\label{eqdefs}
S(\gamma;x,\lambda) = \frac{N(\gamma)\psi(\gamma;x,\lambda)}
   {\sqrt{w_1w_2}\Delta(\gamma;x)\Delta(\hg;\lambda)}.
\end{equation}
The function $N$ \eqref{eqdefn} is a convenient normalization factor when 
matching $S$ to $R$ in Section \ref{sec8}.
\begin{lemma}\label{lempolesS}
$S(\gamma;x,\lambda)$ is meromorphic on $\mc{O}$ with possible poles at
\begin{equation*}
 \lambda = \pm (\nu - i \hg_k),  \qquad 
 x = \pm (\nu - i \gamma_k), \qquad
  i\gamma_0 + i \gamma_l = -\nu  -iw
\end{equation*}
 for $\nu \in \Lambda_+$,  $k = 0,1,2,3$, and $l=1,2,3$.
\end{lemma}
\begin{proof}
Using \eqref{eqdefPsi} and \eqref{eqeandg} we can express $S$ as
\[
S(\gamma;x,\lambda) = 
\frac{\Psi(\gamma;x,\lambda) N(\gamma) }{
\prod_{k=0}^3 E(\pm x + i \gamma_k)E(\pm \lambda + 
     i \hat \gamma_k)} .
\]
From this expression we can easily read off that the possible pole
hyperplanes are as stated in the lemma (they have to be either poles
of $N(\gamma)$ or zeros of one of the $E$-functions in the denominator).
\end{proof}

\begin{theo}\label{lemaweq}
The function $S(\gamma;x,\lambda)$
is a simultaneous eigenfunction of the two Askey-Wilson type
second order difference operators $\mc{L}_\gamma^x$ and $\tc{L}_\gamma^x$ 
\textup{(}see \eqref{eqdefl}\textup{)} with eigenvalues 
$v(\lambda;w_1,w_2,\gamma)$ and $v(\lambda;w_2,w_1,\gamma)$
respectively, where $v$ is defined by \eqref{eqdefv}.
\end{theo}
\begin{proof} 
Note that $\Delta$ satisfies the first order 
difference equation
\begin{equation*}
\Delta(x+iw_1/2) = \frac{\cosh(\frac{\pi}{w_2}(x+i\gamma_2))
                         \cosh(\frac{\pi}{w_2}(x+i\gamma_3))}
                        {\cosh(\frac{\pi}{w_2}(x-i\gamma_0))
                         \cosh(\frac{\pi}{w_2}(x-i\gamma_1))}
                      \Delta(x-iw_1/2). 
\end{equation*}
The desired eigenvalue equation 
\eqref{eqdefl} for $\mc{L}_\gamma^x$ now follows 
immediately from Lemma \ref{lemsecdifeq}. 

To prove the result  for the
operator $\tc{L}_{\gamma}^x$  we note that
$S$ is symmetric in $w_1$ and $w_2$, while interchanging $w_1$ and
$w_2$ transforms $\mathcal{L}$ to $\tilde{\mathcal{L}}$. 
We could also prove the second difference equation
 by repeating the argument for the first difference equation 
using the component $\uqt$ of the modular double.
\end{proof}

We continue the analysis of the eigenfunction $S$ by proving its duality 
in the geometric variable $x$ and the spectral variable $\lambda$,
similar to the duality \eqref{eqrdual} for Ruijsenaars' hypergeometric 
function $R$. The duality transformation $\gamma \to \hg$ of the parameters 
(see \eqref{eqdefhatg}) is equivalent to interchanging $\rho$ and 
$\upsilon$ under the parameter correspondence \eqref{eqmy2ruijs}: 
$(\rho,\sigma,\tau, \upsilon) \to (\upsilon,\sigma,\tau,\rho)$.

\begin{theo}[Duality]\label{lemsdual}
We have
\[
S(\gamma;x,\lambda) = S(\hg;\lambda,x)
\]
as meromorphic functions on $\mc{O}$.
\end{theo}
\begin{proof}
Assume that $w_1,w_2>0$ and $w/2 > \xi + |\Im(x)| + | \Im(\lambda)|$, where
$\xi$ is as in \eqref{eqdefxi}. Note that these restrictions on the parameters
are invariant under the exchange $(x,\gamma) \leftrightarrow (\lambda,\hg)$.
Then we 
can use the integral representation \eqref{eqintpsi} for both
$\psi(\gamma;x,\lambda)$ and $\psi(\hg;\lambda,x)$ to compute 
\begin{equation*}
\begin{split}
\psi(\gamma;x,\lambda) &= 
\int_{\mbb{R}} \frac{G(z+x+\lambda/2-iw/2\pm i\tau)G(z-\lambda/2 -iw/2 \pm 
 i\sigma)}{G(z+x-\lambda/2 + iw/2 \pm i\rho)G(z+\lambda/2+ iw/2 \pm i\upsilon)}
  dz \\
 &= \int_{\mbb{R}} \frac{G(z+x/2 + \lambda - iw/2 \pm i\tau) 
    G(z-x/2 -iw/2 \pm i\sigma)}{G(z+x/2 +iw/2\pm i\rho)
  G(z-x/2 + \lambda + iw/2 \pm i\upsilon)} dz \\
  &= \psi(\hg;\lambda,x),
\end{split}
\end{equation*}
where we used the change of integration variable $z \to z + (\lambda-x)/2$
and a contour shift in the second
equality. This contour shift is allowed since the integrand converges to zero
exponentially
at $\pm \infty$, and the conditions on the parameters ensure that there are no
poles picked up by shifting the contour back to $\mbb{R}$. 

Since $\Psi$ (see \eqref{eqdefPsi}) is entire on $\mc{O}$, it follows that
$\psi(\gamma;x,\lambda) = \psi(\hg;\lambda,x)$ holds as identity between
meromorphic functions on $\mc{O}$.
The desired duality for $S$ now follows from $N(\gamma) = N(\hat \gamma)$ and $\hat{\hat \gamma} = \gamma $.
\end{proof}

\begin{corol}\label{thsawdiff}
The function $S(\gamma;x,\lambda)$ is a simultaneous eigenfunction of the
 Askey-Wilson second order difference operators
$\mc{L}_\gamma^x$, $\tc{L}_\gamma^x$, $\mc{L}_{\hg}^\lambda$, and 
$\tc{L}_{\hg}^\lambda$ with eigenvalues
$v(\lambda;w_1,w_2,\gamma)$, $v(\lambda;w_2,w_1,\gamma)$, 
$v(x;w_1,w_2,\hat \gamma)$, and $v(x;w_2,w_1,\hat \gamma)$
respectively.
\end{corol}
\begin{proof}
The fact that $S$ is an eigenfunction of 
$\mathcal{L}_\gamma^x$ and $\tc{L}_\gamma^x$ 
was proved in Theorem \ref{lemaweq}. The proof for the other two
difference operators follows from this fact and duality
(Theorem \ref{lemsdual}).
\end{proof}

It is immediately clear from the integral 
representation \eqref{eqintpsi} that $\psi$ is invariant under 
sign flips of the
parameters $\rho$, $\sigma$, $\tau$, and $\upsilon$. This leads to the 
following symmetries for $S$.
\begin{lemma}\label{lemSsymmetry}
Let $W_n$ be the Weyl group of type $D_n$, which acts on $n$-tuples by
permutations and even numbers of sign changes. 
Let $V=W_2 \times W_2 \subset W_4$ be the Weyl group of type 
$D_2\times D_2$, where the first
\textup{(}respectively second\textup{)} component acts on the parameters 
$(\gamma_0,\gamma_1)$ \textup{(}respectively $(\gamma_2,\gamma_3)$\textup{)}
 of the four-tuple $(\gamma_0,\gamma_1,\gamma_2,\gamma_3)$.
For an element $v\in V$ we have
\begin{equation*}
\frac{S(\gamma;x,\lambda)}{c(\gamma;x)c(\hat \gamma;\lambda) N(\gamma)}
 = \frac{S(v(\gamma);x,\lambda)}
   {c(v(\gamma);x)c(\widehat{v(\gamma)};\lambda)N(v(\gamma))}
\end{equation*}
as meromorphic functions on $\mc{O}$.
\end{lemma} 
\begin{proof}
Note that the action of $V \simeq \mbb{Z}_2^{\times 4}$ on the parameters 
$(\rho,\sigma,\tau,\upsilon)$ is 
by sign flips of $\rho$, $\sigma$, $\tau$, and $\upsilon$.
Under the conditions $\zeta>0$ and $|\Im(x)|<\zeta$ it follows from the 
integral 
representation \eqref{eqintpsi} of $\psi$ that $\psi$ is invariant
under the action of $V$ on $\gamma$ (note that the parameter restrictions
are $V$-invariant).

Observe that 
\begin{equation*}
c(\gamma;x)\Delta(\gamma;x) = \frac{G(x\pm i\gamma_2)G(x\pm i \gamma_3)}
     {G(2x+iw)}
\end{equation*}
is also $V$-invariant. Since the
action of $V$ commutes with taking dual parameters (which is obvious in the 
parameters $\rho,\sigma,\tau,\upsilon$, 
since $V$ acts by flipping signs while taking dual parameters amounts to 
interchanging $\rho$ and $\upsilon$) we have a similar result for
$c(\hg;\lambda)\Delta(\hg;\lambda)$. Combining these results and using
\eqref{eqdefs} now yields the desired symmetry of 
$S$ for the restricted parameter set. These extra conditions on the 
parameters can be removed by
analytic continuation (compare with the proof of Theorem \ref{lemsdual}).
\end{proof}
\begin{remark}
The symmetries described in Lemma \ref{lemSsymmetry}
should be compared to the $D_4$ symmetry \eqref{eqd4sym} 
of $R$. Note that for $R$ only an $S_3\subset W_4$ symmetry 
holds trivially from its integral representation \eqref{eqdefR}, 
where $S_3$ acts by permuting $\gamma_1$, $\gamma_2$, and $\gamma_3$.
\end{remark}

Let us now consider the asymptotics of $S$, compare with the asymptotics
\eqref{eqRasymp} of $R$.
\begin{lemma}\label{lemanalyticS}
Let $w_1,w_2\in \mbb{R}_{>0}$, $\gamma \in \mbb{C}^4$, and $\lambda \in
\mbb{C}\setminus \mbb{R}$ such that $\zeta >0$, where $\zeta$ is given by
\eqref{eqdefzeta}.
Then
\[
S(\gamma;x,\lambda) = \mathcal{O}(e^{\alpha ( |\Im(\lambda)| 
     - \Re(\hat \gamma_0)-w)|\Re(x)| })
\]
for $\Re(x) \to \pm \infty$,
uniformly for $\Im(x)$ in compact subsets of $(-\zeta,\zeta)$,
where $\alpha=2\pi/w_1w_2$.
\end{lemma}
\begin{proof}
Under the parameter restrictions as stated in the lemma, 
 $S$ does not have $x$-independent poles (see Lemma
\ref{lempolesS}) and the integral representation \eqref{eqintpsi} for $\psi$ 
holds.

In view of \eqref{eqdefs} and the asymptotics
\begin{equation}\label{eqasympdelta}
\frac{1}{\Delta(\gamma;x)} = \mc{O}(
e^{\mp \alpha \hg_0 x})
\end{equation}
for $\Re(x) \to \pm \infty$,
uniformly for $\Im(x)$ in compacta, it suffices to
prove 
\begin{equation}\label{eqasymppsi}
\psi(\gamma;x,\lambda) = \mc{O}(e^{\alpha (|\Im(\lambda)|-w)|\Re(x)|})
\end{equation}
for $\Re(x)\to \pm \infty$, uniformly for $\Im(x)$ in compacta of 
$(-\zeta,\zeta)$.
 The asymptotic formula
\eqref{eqasympdelta} follows  directly
from the estimates \eqref{eqlimitgfrac} and \eqref{eqpropg}
 for the hyperbolic gamma function.

Note that it suffices to prove \eqref{eqasymppsi}
for $\Re(x) \to \infty$
since 
\begin{equation}\label{eqpsisym}
\psi(\gamma;x,\lambda) = \psi(\check \gamma;-x,-\lambda)
\end{equation}
where $\check \gamma = (\sigma,\rho,\upsilon,\tau)$ (in the $\gamma_\mu$ 
notation, $\check \gamma = (-\gamma_0,\gamma_1,\gamma_2,-\gamma_3)$). 
Equation \eqref{eqpsisym} follows by the change of  integration variable
$z \to -z$ in \eqref{eqintpsi} and a subsequent contour shift.

To prove \eqref{eqasymppsi} for $\Re(x) \to \infty$ we 
consider the integral representation \eqref{eqintpsi} of $\psi$. We define
\begin{equation}\label{eqdefepsilon}
\epsilon = \max(w_1,w_2) + \frac12 |\Re(\lambda)| + \max(|\Im(\rho)|,
  |\Im(\sigma)|,|\Im(\tau)|,|\Im(\upsilon)|)
\end{equation}
and we consider the division of $\mbb{R}$ in five intervals
\begin{equation}\label{eqdefintervals}
\begin{split}
I_1 &= (-\infty, -\Re(x) - \epsilon), \qquad
I_2 = (-\Re(x)-\epsilon, - \Re(x) + \epsilon),  \\
I_3 &= (-\Re(x)+\epsilon,-\epsilon), \qquad
I_4 = (-\epsilon,\epsilon), \qquad
I_5 = (\epsilon,\infty),
\end{split}
\end{equation}
for $\Re(x)>2\epsilon$.
We write the integral \eqref{eqintpsi} 
defining $\psi$ as the sum of five integrals over $I_j$ ($j=1,2,\ldots,5$) 
and we bound the integral over each $I_j$ seperately.
The intervals are chosen in such a way that the estimates
\eqref{eqlimitgfrac} and \eqref{eqpropg}
 for the hyperbolic gamma function can be used to 
bound the integrand over the intervals $I_1$, $I_3$, and $I_5$. To estimate
the integrals over the remaining intervals $I_2$ and $I_4$ we use the
fact that their lenghts are finite and independent of $\Re(x)$.
For each interval $I_j$ we show that the integral over $I_j$
is $ \mathcal{O}(e^{\alpha(|\Im(\lambda)|- w) |\Re(x)|})$ 
as $\Re(x) \to \infty$, uniformly for $\Im(x)$ in compact subsets of
$(-\zeta,\zeta)$.  As a consequence $\psi$ is also of this order.
Details are given in Appendix \ref{appa2}.
\end{proof}


\section{Reduction to Askey-Wilson polynomials}\label{sec7}
Using an indirect method, Ruijsenaars \cite[Theorem 3.2]{Ruijs2} 
proved that $R$ reduces to the Askey-Wilson polynomials
\cite{Askey} when the spectral parameter is specialized to certain specific
discrete values. We now show by a direct calculation that $S$ \eqref{eqdefs}
reduces to the Askey-Wilson polynomials for the same discrete spectral values.

Let us first introduce some standard notations for basic hypergeometric 
series, see \cite{Gasper}. For $q \in \mbb{C}$ we write
\begin{align*}
(a;q)_n = \prod_{k=0}^{n-1} (1-aq^k), \\
(a_1,a_2,\ldots,a_k;q)_n = \prod_{j=1}^k (a_j;q)_n.
\end{align*}
The $q$-hypergeometric series is defined by
\[
\ _{s+1}\phi_s\left[ \begin{array}{c}
a_1,\ldots,a_{s+1} \\ b_1,\ldots,b_{s} \end{array} ; q,z \right]
= \sum_{k=0}^\infty \frac{(a_1,\ldots,a_{s+1};q)_k}
{(b_1,\ldots,b_s,q;q)_k} z^k
\]
provided that either $|q|<1$ or that the series terminates.
The Askey-Wilson polynomials \cite{Askey} are defined as 
\[
r_n(x;a,b,c,d~|~q) = 
\ _4\phi_3\left[ \begin{array}{c}
q^{-n},abcdq^{n-1},ae^{2\pi x/w_2},ae^{-2\pi x/w_2} \\
ab,ac,ad \end{array}; q,q \right].
\]
Notice that the $q^{-n}$ term in the above expression causes
the series to terminate. This implies that $r_n$
is a polynomial of degree $n$ in $\cosh(2\pi x/w_2)$. 
Finally, if we use the parameter correspondence 
\begin{equation}\label{eqparcor}
a  = - e^{2\pi i\gamma_0/w_2} q, \quad 
b  = - e^{2\pi i\gamma_1/w_2} q, \quad 
c  = - e^{2\pi i\gamma_2/w_2} q, \quad 
d  = - e^{2\pi i\gamma_3/w_2} q,
\end{equation}
and if we define
\begin{equation}\label{eqdeflambdan}
\lambda_n = iw+i\hg_0 + in w_1,
\end{equation}
then the Askey-Wilson polynomials satisfy the Askey-Wilson second order
difference equation
\begin{equation*}
\mathcal{L}_\gamma^x r_n(x;a,b,c,d~|~q^2) = 
 v(\lambda_n;w_1,w_2,\gamma)
r_n(x;a,b,c,d ~|~q^2).
\end{equation*}
Here we use the Askey-Wilson operator $\mathcal{L}_\gamma^x$ 
\eqref{eqdefl} and eigenvalue $v$ \eqref{eqdefv}.

Ruijsenaars has shown in \cite{Ruijs2} by an indirect method that 
\begin{equation}\label{eqawr}
R(w_1,w_2,\gamma;x,\lambda_n) = r_n(x;a,b,c,d~|~q^2)
\end{equation}
for $n\in \mbb{Z}_{\geq 0}$,
under the parameter correspondence \eqref{eqparcor}. Similarly we have
\begin{theo}\label{theoaw}
Under the parameter correspondence \eqref{eqparcor} we have
\[
S(w_1,w_2,\gamma;x,\lambda_n) = r_n(x;a,b,c,d~|~q^2)
\]
for $n \in \mbb{Z}_{\geq 0}$.
\end{theo}
\begin{proof}
Without loss of generality we assume that the parameters
$w_1,w_2,\gamma,x$ are generic.

For generic $\lambda$ we can express $\psi$ as an integral 
\begin{equation}\label{eqjasper1}
\psi(\gamma;x,\lambda) = \int_{\mc{C}} I(\gamma;x,\lambda,z)dz
\end{equation}
with $I(z)=I(\gamma;x,\lambda,z)$ given by
\[
I(z)= \frac{G(z+x+\lambda/2-iw/2\pm i \tau)G(z-\lambda/2-iw/2\pm i \sigma)}
     {G(z+x-\lambda/2+iw/2\pm i \rho)G(z+\lambda/2+iw/2\pm i \upsilon)}
\]
and with  contour $\mc{C}$ a deformation of $\mbb{R}$ 
seperating the upward pole sequences of $I$ from the downward pole sequences
of $I$.
When $\lambda \to \lambda_n$, the pole $z_k := \lambda/2 - iw/2-i\sigma-ikw_1$
from a downward pole sequence of $I$ will collide with the pole
$-\lambda/2 +iw/2-i\upsilon + i(n-k)w_1$ from an upward
pole sequence of $I$ for $0\leq k \leq n$. In order to 
compute the limit $\lambda \to \lambda_n$ in \eqref{eqjasper1}, we 
therefore first  
shift the contour $\mc{C}$ over the poles at $z_k$ ($0\leq k\leq n$) while
picking up poles. In the resulting integral the colliding poles are on the
same side of the integration contour, hence the limit $\lambda \to \lambda_n$
can be taken.

To calculate the residues of $I$ at $z_k$
we first remark that $k$ consecutive applications of 
the difference equation \eqref{eqdiffG} yield
\begin{equation*}
\frac{G(z)}{G(z-ikw_1)} = 
  e^{k\pi z/w_2}q^{-k^2/2} (-e^{-2\pi z/w_2}q;q^2)_k.
\end{equation*}
Using this equation we can write 
\begin{align*}
I(z)&=  
\frac{G(z +ikw_1+x+\lambda/2-iw/2\pm i \tau)G(z+ikw_1-\lambda/2-iw/2\pm i \sigma)}
     {G(z +ikw_1 +x-\lambda/2+iw/2\pm i \rho)G(z+ikw_1+\lambda/2+iw/2\pm i \upsilon)}  \\
& \qquad \times 
q^{2k}\frac{(-e^{-\frac{2\pi}{w_2}(z+ikw_1+x-\lambda/2+iw/2 \pm i \rho)}q,
-e^{-\frac{2\pi}{w_2}(z+ikw_1+\lambda/2+iw/2\pm i \upsilon)}q
;q^2)_k}
{(-e^{-\frac{2\pi}{w_2}(z+ikw_1+x+\lambda/2-iw/2 \pm i \tau)}q,
-e^{-\frac{2\pi}{w_2}(z+ikw_1-\lambda/2-iw/2\pm i \sigma)}q,q^2)_k}.
\end{align*}
Using the fact that the residue of the hyperbolic gamma function at
$z=-iw$ equals \eqref{eqresg}, we obtain that the residue $Res_k$ of $I$ at 
$z_k$ equals
\begin{align*}
Res_k &= \frac{i\sqrt{w_1w_2}}{2\pi} 
\frac{G(x+\lambda-iw -i\sigma \pm i \tau)G(-iw-2i\sigma)}
{G(x -i\sigma \pm i\rho)G(\lambda -i\sigma \pm i \upsilon)} \\
& \qquad \times q^{2k} 
\frac{(-e^{-\frac{2\pi}{w_2}(x - i\sigma \pm i \rho)}q,
  -e^{-\frac{2\pi}{w_2}(\lambda-i\sigma \pm i \upsilon)}q;q^2)_k}
{(e^{-\frac{2\pi}{w_2}(x+\lambda-i\sigma \pm i \tau)}q^2,
e^{\frac{2\pi}{w_2} 2i\sigma}q^2,q^2;q^2)_k}.
\end{align*}
Now we can rewrite $S$ as
\begin{equation*}
S(\gamma;x,\lambda) = \frac{N(\gamma)}{\sqrt{w_1w_2}\Delta(\gamma;x)
\Delta(\hg;\lambda)}\left( -2\pi i \sum_{k=0}^n Res_k + 
\int_{\mc{C}^\prime} I(z) dz \right),
\end{equation*}
where the contour $\mc{C}^\prime$ is chosen in such a way 
that all  upward pole sequences and the poles $z_k$ ($0\leq k\leq n$) are
above $\mc{C}^\prime$, while all poles in downward pole sequences except
$z_k$ ($0\leq k \leq n$) are below $\mc{C}^\prime$.
In this expression the integral $\int_{\mc{C}^{\prime}} I(z)dz$
has an analytic extension to $\lambda = \lambda_n$. 
Furthermore $S(\gamma;x,\lambda)$ is analytic at $\lambda = \lambda_n$, while
$\Delta(\hg;\lambda)$ and $Res_k$ ($0\leq k \leq n$) have simple poles
at $\lambda = \lambda_n$. Hence we obtain
\begin{align*}
S(\gamma;x,\lambda_n) &= \lim_{\lambda \to \lambda_n}
-\frac{2\pi i N(\gamma)}{\sqrt{w_1w_2}\Delta(\gamma;x)\Delta(\hg;\lambda)}
\sum_{k=0}^n Res_k \\
& = 
e^{\frac{2n\pi}{w_2}(x-iw-i\gamma_0)} 
\frac{(e^{-\frac{2\pi}{w_2}(x-iw+\gamma_{2/3})}q^{-2n};q^2)_n}
{(e^{-\frac{2\pi}{w_2}(i\gamma_0 + i \gamma_{2/3})}q^{-2n};q^2)_n} \\
& \quad \times \ _4\phi_3\left[ \begin{array}{c}
q^{-2n}, e^{-\frac{2\pi}{w_2}(x-iw-i\gamma_{0/1})}, 
e^{-\frac{2\pi i}{w_2}(\gamma_2+\gamma_3)}q^{-2n} \\
e^{-\frac{2\pi}{w_2}(x-iw + i \gamma_{2/3})}q^{-2n},
e^{\frac{2\pi i}{w_2}(\gamma_0 + \gamma_1)}q^2 \end{array} 
; q^2,q^2 \right],
\end{align*}
where the notation $\gamma_{0/1}$ (respectively $\gamma_{2/3}$)
means that there are two terms, one with $\gamma_0$ and 
another with $\gamma_1$ (respectively, $\gamma_2$ and $\gamma_3$).
Inserting the parameter correspondence \eqref{eqparcor} we obtain
\begin{align*}
S(\gamma;x,\lambda_n) &= 
 e^{2\pi n x/w_2}a^{-n} 
\frac{(e^{-2\pi x/w_2}c^{-1}q^{-2n+2},e^{-2\pi x/w_2}d^{-1}q^{-2n+2};q^2)_n}{(a^{-1}c^{-1}q^{-2n+2},a^{-1}d^{-1}q^{-2n+2};q^2)_n} \\
& \qquad \times\ _4\phi_3\left[ \begin{array}{c}
q^{-2n},  e^{-2\pi x/w_2}a, e^{-2\pi x/w_2}b,c^{-1}d^{-1}q^{-2n+2} \\
 e^{-2\pi x/w_2}c^{-1}q^{-2n+2}, e^{-2\pi x/w_2}d^{-1}q^{-2n+2},ab \end{array}
;q^2,q^2 \right].
\end{align*}
Using Sears' transformation \cite[(III.15)]{Gasper}
of a terminating balanced $_4\phi_3$ series  with parameters 
specialized to $a = a e^{-2\pi x/w_2}$, $b=b e^{-2\pi x/w_2}$, $c=c^{-1}d^{-1}q^{-2n+2}$, 
$d=ab$, $e=e^{-2\pi x/w_2}c^{-1}q^{-2n+2}$, and $f=e^{-2\pi x/w_2}d^{-1}q^{-2n+2}$ now yields the desired result.
\end{proof}


\section{Equality to Ruijsenaars' hypergeometric function}\label{sec8}

We have already seen in previous sections that Ruijsenaars' hypergeometric
function $R$ and the renormalized formal matrix coefficient $S$
have several properties in common. They satisfy the same Askey-Wilson
second order difference equations, they have the same duality property, 
they specialize in the same way to the Askey-Wilson polynomials and 
their possible pole locations coincide. These common properties suffice
to show that $R$ and $S$ are equal.
\begin{theo}\label{thgmain}
We have
\begin{equation}\label{eqthgmain}
R(w_1,w_2,\gamma;x,\lambda) = S(w_1,w_2,\gamma;x,\lambda).
\end{equation}
\end{theo}

This theorem is equivalent to the following identity between hyperbolic 
integrals.
\begin{corol}
For $w_1,w_2, \Re(\gamma_j)>0$ and $|x|,|\lambda|,|\gamma_j| < w/6$ we have
\begin{equation*}
\begin{split}
\int_{\mbb{R}}&
  \frac{G(z+x+\lambda/2 -iw/2 \pm i(\gamma_3-\gamma_2)/2) 
        G(z-\lambda/2 -iw/2 \pm i (\gamma_0+\gamma_1)/2)}
   {G(z+x-\lambda/2 +iw/2 \pm i (\gamma_0-\gamma_1)/2)
 G(z+\lambda/2 +iw/2 \pm i (\gamma_2+\gamma_3)/2)} dz \\
&= \frac{ G(x+i\gamma_2)G(x+i\gamma_3)G(\lambda +i\hg_2)G(\lambda+i\hg_3)}
  {G(x+i\gamma_0)G(x-i\gamma_1)G(\lambda+i\hg_0)G(\lambda - i \hg_1) } \\
 &\qquad \qquad  \qquad \times \int_{\mc{C}}
  \frac{G(z \pm x + i\gamma_0) G(z \pm \lambda + i \hg_0)}
  {G(z+ iw)
 \prod_{j=1}^3 
            G(z+ i\gamma_0 + i \gamma_j + iw)} dz ,
\end{split}
\end{equation*}
where the contour $\mc{C}$ is the real line with a downward 
indentation at the origin.
\end{corol}
\begin{proof}
The proof consists of inserting the integral representations 
of $R$ and $S$ in \eqref{eqthgmain}. See \eqref{eqdefR} for the integral
representation of $R$, and \eqref{eqintpsi}, \eqref{eqdefs} for the
integral representation of $S$.
\end{proof}

In order to prove Theorem \ref{thgmain} we first consider 
the Casorati-determinant of $S$ and $R$ in the $iw_1$ direction.
\begin{lemma}\label{lemCas}
The Casorati-determinant
\begin{multline*}
\delta(\gamma;z,\lambda) =\\  S(\gamma;z+iw_1/2,\lambda)R(\gamma;z-iw_1/2,\lambda)-
S(\gamma;z-iw_1/2,\lambda)R(\gamma;z+iw_1/2,\lambda)
\end{multline*}
of $S$ and $R$ in the $iw_1$ direction is identically zero.
\end{lemma}
\begin{proof}
We suppress the $\lambda$ and $\gamma$ dependence of $\delta(z)$ whenever this
does not cause confusion. We prove the lemma for generic parameters
$w_1,w_2\in \mbb{R}_{>0}$, $\gamma \in \mbb{R}^4$, and 
$\lambda \in U \setminus \mbb{R}$, under the condition
$w_2> 2 \xi +2 |\Im(\lambda)| +3w_1$, where $U$ is an open subset such that
the asymptotics \eqref{eqRasymp} of $R$ hold for $\lambda \in U$.

A simple calculation involving the Askey-Wilson difference
equations satisfied by $R$ and $S$ 
(see \eqref{eqrdiff} and Theorem \ref{thsawdiff} respectively) shows that 
\[
\delta(z+iw_1/2) =
  \frac{A(\gamma;-z)}{A(\gamma;z)} \delta(z-iw_1/2),
\]
where $A$ is defined by \eqref{eqdefar}.
Since the function 
\begin{equation*}
T(z) = \sinh(2\pi z/w_2) \prod_{j=0}^3 
    \frac{G(z-i\gamma_j-iw_1/2)}{G(z+i\gamma_j+iw_1/2)}
\end{equation*}
 satisfies the same difference equation, we conclude that
\[
m(z) = \frac{\delta(z)}{T(z)}
\]
is an $iw_1$-periodic function. 

We now show that $m(z)$ is an
entire function in $z$.
Let us look at the possible poles of the Casorati-determinant
$\delta(z)$. By Lemma \ref{lemanalyticS} the possible poles of $S$ are 
located at
\[
\pm (\Lambda_+ - i \gamma_j),  \qquad (j = 0,1,2,3).
\]
From \eqref{eqpolesR} 
the possible poles of $R$ are located at the same points. Hence
 $\delta(z)$ can only have poles at
\[
  \pm(\Lambda_+  - i \gamma_j) \pm iw_1/2 , \qquad (j = 0,1,2,3)
\]
Here all sign combinations are possible.
Furthermore, using the pole and zero locations \eqref{eqdeflambdaplus} of the
hyperbolic gamma function,
we can easily see that the possible zeros of $T(z)$ are located at
\[
 \pm (\Lambda_+ + i\gamma_j+iw_1/2), \qquad   riw_2, \qquad
 (j =0,1,2,3;\; r\in \mbb{Z}).
\]
By the assumption that the parameters are generic, we conclude that $m$ has
no pole sequences of the form $p+ikw_1$ ($k\in \mbb{Z}$).
By the $iw_1$-periodicity of $m$ it now follows that $m$ cannot have any 
poles. 

In the limit  $\Re(z) \to \infty$ we have
\[
\frac{1}{T(z)} = \mc{O}( e^{\alpha (\hg_0+w_1)z})
\]
uniformly for $\Im(z)$ in compacta,
in view of the estimates \eqref{eqlimitgfrac} and \eqref{eqpropg}
for the hyperbolic gamma function. Here $\alpha = 2\pi/w_1w_2$ as before. 

Furthermore, using the asymptotics for $S$ (see  Lemma \ref{lemanalyticS}) and 
for $R$ (see \eqref{eqRasymp}) we have for $\Re(z) \to \infty$ 
\[
\delta (z) = 
\mc{O}(e^{2\alpha(|\Im(\lambda)|+|\hg_0| - w)|\Re(z)|})
\]
uniformly for $\Im(z)$ in compact subsets of $(-\zeta+w_1/2,\zeta-w_1/2)$.
Observe that the interval $(-\zeta+w_1/2,\zeta-w_1/2)$ is nonempty due
to the conditions on the parameters. 

Combining these two asymptotic estimates we obtain 
\begin{equation}\label{eqasympm}
m(z) = \frac{\delta(z)}{T(z)} = \mc{O}(e^{\alpha(2|\Im(\lambda)| -w_2)|\Re(z)|}) \to 0
\end{equation}
for $\Re(z) \to \infty$, uniformly for $\Im(z)$ in compacta of
$(-\zeta+w_1/2, \zeta - w_1/2)$. 

The asymptotics of $m(z)$ for $\Re(z) \to -\infty$ can be obtained in a similar
way and is also given by \eqref{eqasympm}. 
Combining the asymptotics with the fact that $m(z)$ is analytic and 
$iw_1$-periodic we conclude that $m(z)$ is bounded on $\mbb{C}$ since
$\zeta-w_1/2 > w_1/2$. 

For these parameters we conclude by Liouville's theorem
that $m(z)$ is constant. In fact, by the asymptotic expansion
\eqref{eqasympm}, $m$ is identically zero. 
We can now extend this result to all
values of the parameters by analytic continuation, which proves the lemma.
\end{proof}

\begin{proof}[Proof of Theorem \ref{thgmain}]
Consider the quotient
\[
Q(\gamma;x,\lambda) = \frac{R(\gamma;x,\lambda)}{S(\gamma;x,\lambda)}.
\]
By Lemma \ref{lemCas}, $Q$ is an $iw_1$-periodic meromorphic function in $x$.
Since $Q$ is symmetric in $w_1$ and $w_2$ (for both $R$ and $S$ are
invariant under interchanging $w_1$ and $w_2$), $Q$ is 
also $iw_2$-periodic. If we choose $w_1,w_2>0$ such that 
$w_1/w_2 \not \in \mbb{Q}$, then the set
$\{kw_1+lw_2~|~k,l\in \mbb{Z}\}$ is dense on the real line, hence 
$Q(\gamma;x,\lambda)$ is constant as meromorphic function in $x$. Analytic 
continuation
(in $w_1$, $w_2$, and $\gamma$) allows us to extend this result
to all possible values of $w_1$ and $w_2$ in $\mbb{C}_+$ and 
$\gamma \in \mbb{C}^4$. 

By the duality properties of $R$ and $S$ (see
\eqref{eqrdual} and Theorem \ref{lemsdual} respectively), we have
\[
Q(\gamma;x,\lambda) = Q(\hg;\lambda,x).
\]
This implies that $Q$ is also constant as function in $\lambda$.

In particular we have
\[
Q(w_1,w_2,\gamma;x,\lambda) = \frac{S(w_1,w_2,\gamma;x,\lambda_0)}
 {R(w_1,w_2,\gamma;x,\lambda_0)}
\]
with $\lambda_0$ given by \eqref{eqdeflambdan}. By Theorem \ref{theoaw}
we have $S(w_1,w_2,\gamma;x,\lambda_0) \equiv 1$, and by 
\eqref{eqawr} we have $R(w_1,w_2,\gamma;x,\lambda_0)\equiv 1$. 
Hence $Q\equiv 1$, as desired.
\end{proof}

\appendix


\section{Eigenfunction of  $\phl(Y_\rho)$}\label{app1}
In this appendix we give the explicit calculation to rewrite the 
eigenvalue equation $\phl(Y_\rho)f = \mu_{\tau}(\rho) f$ as the
 first order difference equation \eqref{eqfode}. 

Using the explicit expression \eqref{eqdefyrho} of 
 $Y_{\rho}$, the eigenvalue equation becomes
\[
iq^{-1/2}\phl(E)f + iq^{-1/2}\phl(FK)f - \frac{\nu_{\rho}}{q-q^{-1}}
(\phl(K-1))f = \frac{\nu_{\rho}-\nu_{\tau}}{q-q^{-1}} f.
\]
By the explicit definition (Lemma \ref{lemrep}) of $\phl$ we obtain
\begin{equation*}
\begin{split}
\frac{i}{q-q^{-1}}e^{2\pi z/w_2}
\left( 
q^{-1/2}e^{\pi \lambda/w_2}f(z) + q^{1/2}e^{-\pi\lambda/w_2}f(z+iw_1)
\right) & \\
-\frac{i}{q-q^{-1}}e^{-2\pi z /w_2} 
\left(
q^{-1/2}e^{\pi\lambda/w_2} f(z+iw_1) + q^{1/2}e^{-\pi\lambda/w_2}f(z)
\right) & \\
-\frac{\nu_\rho}{q-q^{-1}}\left( f(z+iw_1)-f(z)\right) = &
\frac{\nu_{\rho}-\nu_\tau}{q-q^{-1}}f(z).
\end{split}
\end{equation*}
Multiplying by $q-q^{-1}$ and rearranging the terms yields
\begin{equation*}
\begin{split}
& \left(i  e^{2\pi z /w_2}q^{-1/2}e^{\pi \lambda/w_2}
- i e^{-2\pi z/w_2}q^{1/2}e^{-\pi\lambda/w_2} + \nu_\tau \right) f(z) \\
&= \left(
-ie^{2\pi z/w_2} q^{1/2}e^{-\pi\lambda/w_2} +i e^{-2\pi z /w_2}q^{-1/2}
e^{\pi\lambda/w_2} + \nu_\rho \right) f(z+iw_1),
\end{split}
\end{equation*}
which is equivalent to
\begin{equation*}
f(z+iw_1) = 
\frac{ \cosh(\pi i/2 + 2\pi z/w_2  - \pi i w_1/(2w_2)+ \pi \lambda/w_2)
   +\cosh(2 \pi i\tau/w_2)}
{\cosh(-\pi i/2 +2\pi z/w_2 + \pi iw_1/(2w_2) -\pi\lambda/w_2)
   + \cosh(2 \pi i \rho/w_2)} f(z).
\end{equation*}
Replacing the variable $z$ by $z-iw_1/2$ we can now rewrite the latter equation
as
\begin{equation*}
\begin{split}
&f(z+iw_1/2) \\
&= 
\frac{\cosh(\pi i/2 + 2\pi z/w_2 -3\pi i w_1/(2w_2) + \pi \lambda /w_2)
  + \cosh(2 \pi i \tau/w_2)}
{\cosh(-\pi i/2 + 2\pi z /w_2  - \pi i w_1/(2w_2) - \pi \lambda/w_2)
+ \cosh(2\pi i \rho /w_2)} f(z-iw_1/2) \\
&= \frac{\cosh(\frac{\pi}{w_2}
  (z +\lambda/ 2-3iw_1/4 + iw_2/4 \pm i\tau))}{
 \cosh(\frac{\pi}{w_2}(z - \lambda/2 - iw_1/4 - iw_2/4 \pm i\rho))}
 f(z-i w_1/2) \\
&=  \frac{\cosh(\frac{\pi}{w_2}(z +\lambda /2-3iw/2 \pm i\tau))}
  {\cosh(\frac{\pi}{w_2}(z - \lambda/2 - iw/2 + i\rho))} f(z-iw_1/2), 
\end{split}
\end{equation*}
where we used the $i\pi$-antiperiodicity of the hyperbolic
cosine in the last equality. 


\section{The limit behaviour of $\psi$}\label{appa2}

In this appendix we give the details on the calculation of the 
limit behaviour of $\psi$, cf.\ the proof of Lemma \ref{lemanalyticS}. 
Throughout this 
section we assume that $w_1,w_2\in \mbb{R}_{>0}$, 
$\lambda \not\in \mbb{R}$,
$\gamma \in \mbb{C}^4$, and that
$\zeta>0$ (with $\zeta$ given by \eqref{eqdefzeta}).
We prove that 
\begin{equation}\label{eqlimpsi}
\psi(\gamma;x,\lambda) = \mc{O}(e^{\alpha(|\Im(\lambda)|-w)|\Re(x)|})
\end{equation}
for $\Re(x) \to \infty$, uniformly for $\Im(x)$ in compacta of 
$(-\zeta,\zeta)$.
As explained in the proof of Lemma \ref{lemanalyticS}, we prove 
\eqref{eqlimpsi} by splitting $\mbb{R}$ in five intervals and bounding the 
integral representation \eqref{eqintpsi} of $\psi$ over each interval.

\subsection{Preparations}

Let us first define a function $K$ by
\[
K(z,\lambda,a,b) = 
\frac{G(z+\lambda/2-iw/2\pm ia)}
     {G(z-\lambda/2+iw/2 \pm ib)}.
\]
The integral representation \eqref{eqintpsi} for $\psi$
can then be written as
\begin{equation}\label{eqintpsi2}
\psi(\gamma;x,\lambda)  =
\int_{\mbb{R}} K(z+x,\lambda,\tau,\rho)K(z,-\lambda,\sigma,\upsilon)dz.
\end{equation}
The behaviour of $K$ in the limit $z \to \pm \infty$ is controlled by
\begin{equation*}
 K_\pm(z,\lambda,a,b) = e^{\mp i\alpha(z(\lambda-iw)-a^2/2+b^2/2)}.
\end{equation*}
Explicitely, for fixed $a$, $b$, and $\lambda$ we have
\begin{equation}\label{eqlimitk+}
K(z,\lambda,a,b) = K_\pm(z,\lambda,a,b)e^{g(z,\lambda,a,b)}
\end{equation}
for $\pm \Re(z)> \max(w_1,w_2) + |\Re(\lambda)|/2+ \max(|\Im(a)|,|\Im(b)|)$,
where 
\begin{equation}\label{eqcondf}
|g(z,\lambda,a,b)| < C(\Im(z)) e^{-\alpha \min(w_1,w_2)|\Re(z)|/2},
\end{equation}
with $C$ depending continuously on $\Im(z)$,
cf.\ \eqref{eqlimitgfrac} and \eqref{eqpropg}.

\subsection{General estimation scheme}

Let $\epsilon$ and the intervals $I_j$ ($j \in \{1,\ldots,5\}$) be defined
as in \eqref{eqdefepsilon} and \eqref{eqdefintervals}. We only
consider the asymptotics for $\Re(x) \to \infty$. 
Assume that $\Re(x)>2\epsilon$, causing the intervals to form a partition of 
the real line. We write the integral \eqref{eqintpsi2} defining $\psi$ as 
\begin{equation}\label{eqsplitpsi}
\psi(x) = \sum_{j=1}^5 \psi_j(x),
\end{equation}
where
\begin{equation*}
\psi_j(x) = 
   \int_{I_j} K(z+x,\lambda,\tau,\rho)K(z,-\lambda,\sigma,\upsilon)dz
\end{equation*}
for $j\in \{1,2,\ldots,5\}$. 
We bound these intergrals using \eqref{eqlimitk+}
(if one of them is applicable for the interval at hand).

For $j=1$ we have
\begin{equation*}
\begin{split}
\psi_1(x) 
& = \int_{-\infty}^{-\Re(x)-\epsilon} K_-(z+x,\lambda,\tau,\rho) 
   K_-(z,-\lambda,\sigma,\upsilon) e^{g_1(z+x,x)}dz   \\
&= e^{i\alpha \lambda x} e^{-\alpha w \bar x}
e^{i\alpha(\rho^2+\upsilon^2-\tau^2-\sigma^2)/2}
  \int_{-\infty}^{-\epsilon} e^{2\alpha w z +
       g_1(z+i\Im(x),x)}dz   \\
&= \mc{O}(e^{-\alpha (\Im(\lambda)+w) \Re(x)})
\end{split}
\end{equation*}
for $\Re(x) \to \infty$, uniformly for $\Im(x)$ in compacta of 
$(-\zeta,\zeta)$.
Here $g_1(z,x)=g(z,\lambda,\tau,\rho)+g(z-x,-\lambda,\sigma,\upsilon)$
which satisfies an equation like \eqref{eqcondf} for $z<-\epsilon$ 
\[
|g_1(z+i\Im(x),x)| < C e^{-\alpha \min(w_1,w_2)|\Re(z)|/2}.
\]
where the constant $C$ is independent of $\Im(x)$, because 
$\Im(x)$ is bounded. 
In particular, $g_1(z+i\Im(x),x)$ is uniformly bounded for $z\in (-\infty,-\epsilon)$
and $x \in \{z\in \mbb{C} ~|~ \Re(z) \geq 2\epsilon, |\Im(z)|<\zeta\}$.

Likewise we have for $j=5$,
\begin{equation*}
\begin{split}
\psi_5(x) 
& = \int_{\epsilon}^{\infty} K_+(z+x,\lambda,\tau,\rho)
   K_+(z,-\lambda,\sigma,\upsilon)e^{g_5(z,x)} dz 
  \\
&= e^{-\alpha x (w + i\lambda)}
e^{i\alpha(\sigma^2+\tau^2-\rho^2-\upsilon^2)/2}
\int_{\epsilon}^{\infty} e^{-2\alpha w z +g_5(z,x)}dz
\\
&= \mc{O}(e^{\alpha(\Im(\lambda) -w)\Re(x)})
\end{split}
\end{equation*}
for $\Re(x)\to \infty$, uniformly for $\Im(x)$ in compacta of 
$(-\zeta,\zeta)$.
Here $g_5$ is a function which satisfies a bound like 
\eqref{eqcondf} for $z>\epsilon$, cf.\ the previous paragraph.

For $j=3$ we need to be a bit more careful. First observe that 
\begin{equation*}
\begin{split}
\psi_3(x)  &= 
 \int_{-\Re(x)+\epsilon}^{-\epsilon} K_+(z+x,\lambda,\tau,\rho)
  K_-(z,-\lambda,\sigma,\upsilon) e^{g_3(z,x)}dz  \\
&=  e^{-\alpha x(w+i\lambda)}
e^{i\alpha(\tau^2+\upsilon^2-\rho^2-\sigma^2)/2}
\int_{-\Re(x)+\epsilon}^{-\epsilon} e^{-2i\alpha \lambda z + g_3(z,x)} dz,
\end{split}
\end{equation*}
where $g_3 = g(z+x,\lambda,\tau,\rho) + g(z,-\lambda,\sigma,\upsilon)$ is
bounded on $z \in (-\Re(x)+\epsilon,-\epsilon)$ by
$C e^{-\alpha \min(w_1,w_2) \min(-z,z-\Re(x))/2}$, and hence by the constant 
$C$ itself. Therefore  we have
\begin{equation*}
\begin{split}
|\psi_3(x)| 
& \leq  e^{\alpha  (\Im(\lambda)-w)\Re(x) +\alpha \Im(x)\Re(\lambda)}
e^{\alpha \Im(\tau^2+\upsilon^2-\rho^2-\sigma^2)/2}
\int_{-\Re(x)+\epsilon}^{-\epsilon} e^{2\alpha \Im(\lambda) z +C} dz\\
&= \mc{O}(e^{\alpha(|\Im(\lambda)|- w)\Re(x)})
\end{split}
\end{equation*}
for $\Re(x)\to \infty$, uniformly for $\Im(x)$ in compacta of 
$(-\zeta,\zeta)$.
Here we get the final approximation by evaluating the integral and using  
that $\Im(\lambda) \neq 0$.

For $j=4$ we cannot use \eqref{eqlimitk+} for the entire integrand.
However we still have
\begin{equation*}
\begin{split}
\psi_4(x)
 &= \int_{-\epsilon}^{\epsilon} K_+(z+x,\lambda,\tau,\rho)
  K(z,-\lambda,\sigma,\upsilon) e^{g_4(z,x)} dz  \\
 &= e^{i\alpha x(iw-\lambda)}    
  \int_{-\epsilon}^{\epsilon} 
 K_+(z,\lambda,\tau,\rho)
  K(z,-\lambda,\sigma,\upsilon)e^{g_4(z,x)} dz  \\
&= \mc{O}(e^{\alpha(\Im(\lambda)-w)\Re(x)})
\end{split}
\end{equation*}
for $\Re(x)\to \infty$, uniformly for $\Im(x)$ in compacta of
$(-\zeta,\zeta)$.
Here $g_4$ is a function satisfying the bound
$g_4(z,x) < C e^{-\alpha \min(w_1,w_2) (\Re(x)-\epsilon)/2} \leq C$, for
$z\in [-\epsilon,\epsilon]$ and $\Re(x)>2\epsilon$.

Finally for $j=2$ we have in a similar way
\begin{equation*}
\begin{split}
\psi_2(x) 
&= \int_{-\epsilon}^{\epsilon} K(z+i\Im(x),\lambda,\tau,\rho)
  K_-(z-\Re(x),-\lambda,\sigma,\upsilon)e^{g_2(z,x)}dz  \\
&= e^{-i\alpha \Re(x)(-\lambda-iw)} 
  \int_{-\epsilon}^{\epsilon} K(z+\Im(x),\lambda,\tau,\rho)
  K_-(z,-\lambda,\sigma,\upsilon)e^{g_2(z,x)}dz  \\
&= \mc{O}(e^{-\alpha (\Im(\lambda)+w) x})
\end{split}
\end{equation*}
for $\Re(x) \to \infty$, uniformly for $\Im(x)$ in compacta of
$(-\zeta,\zeta)$, where $g_2$ is a bounded function, cf.\ the previous 
paragraph.

By \eqref{eqsplitpsi} we conclude that the asymptotics \eqref{eqlimpsi} for
$\psi$ holds, as desired.

\end{document}